\documentclass[a4paper,11pt,twoside]{amsart}

\usepackage{amsmath,amssymb,amsfonts,amsthm,mathtools}
\usepackage[english]{babel}
\usepackage[utf8]{inputenc}
\usepackage{url}
\usepackage[pdftex,citecolor=green,linkcolor=red]{hyperref}
\usepackage[a4paper,inner=3cm,outer=3cm,top=4cm,bottom=4cm,pdftex]{geometry}
\usepackage{fancyhdr}
\usepackage{comment}

\usepackage{todonotes}
\newtheorem{theorem}{Theorem}[section]
\newtheorem{corollary}[theorem]{Corollary}
\newtheorem{lemma}[theorem]{Lemma}
\newtheorem{proposition}[theorem]{Proposition}

\numberwithin{equation}{section}

\newcommand{\U}{\mathbb U}
\newcommand{\N}{\mathbb N}
\newcommand{\C}{\mathbb C}
\newcommand{\R}{\mathbb R}
\newcommand{\D}{\mathbb D}
\newcommand{\1}{\mathbf 1}
\newcommand{\re}{\operatorname{Re}}

\newcommand{\eps}{\varepsilon}
\newcommand{\logstar}{\log_\ast}

\newcommand{\abs}[1]{\left\lvert {#1} \right \rvert}

\title{Inverse Littlewood problem for multiplicative functions}
\author{P\'ea Bazin \and Oleksiy Klurman}
\date{June 6, 2026}

\begin{document}

\providecommand{\de}{\,\mathrm{d}}
\providecommand{\conj}[1]{\overline{#1}}
\providecommand{\bra}[1]{\left( {#1} \right)}
\providecommand{\Q}{\mathbb{Q}}
\providecommand{\Ind}{\text{\fontencoding{U}\fontfamily{bbold}\fontseries{m}\fontshape{n}\selectfont 1}}

\begin{abstract}
We formulate and solve the asymptotic inverse Littlewood problem
in the class of multiplicative functions. That is, let $f:\N\to\U\cup\{0\}$ be a bounded completely multiplicative function satisfying
$\sum_{p:f(p)=0}1/p<\infty$. We show that
\[
 \int_0^1 \left|\sum_{n\le x}f(n)e(\alpha n)\right|\,d\alpha\ll\log x
\]
for all $x\ge3$ if and only if $f(n)=\chi(n)n^{it}$ for some Dirichlet character
$\chi$ and some $t\in\R$.\\
More generally, we show that there is $\eps_0>0$, depending only on the $\ell^2$-density of $f$ (that is, on $\prod_{p:f(p)=0}(1-1/p)$), such that under the weaker hypothesis
\[
 \int_0^1 \left|\sum_{n\le x}f(n)e(\alpha n)\right|\,d\alpha
 \ll x^{\eps_0}
\]
for all $x\ge3$, there are a Dirichlet character $\chi$ and
$t\in\R$ such that
\[
 \sum_p\frac{\abs{f(p)-\chi(p)p^{it}}}{p}<\infty.
\]
We give evidence for the sharpness of such classification, by constructing  non-pretentious examples with
$L^1$ norm $o(\sqrt{x})$ at every sufficiently large scale, but with ``wild" behavior at the intermediate scales.
\end{abstract}

\maketitle

\section{Introduction}

Throughout the paper, let $\U$ denotes the unit circle, and we write
\[
 e(u):=e^{2\pi i u}.
\]
Littlewood conjectured that there is an absolute constant $c>0$ such that every finite set $S \subset \mathbb{N}$ satisfies
\begin{equation} \label{eq:little}
\int_{0}^{1} \Big | \sum_{n \in S} e^{2\pi i n \alpha} \Big | d \alpha > c \log |S|.
\end{equation}
He further conjectured that, among all $S\subset\mathbb{N}$ with $|S|=N$, the integral in \eqref{eq:little} is smallest for $S=[1,N]$. The first conjecture was proved independently by Konyagin~\cite{Konyagin} and by McGehee, Pigno and Smith~\cite{MPS}, while the second one remains wide open. More generally, the results of \cite{MPS} give the same lower bound for arbitrary coefficients $a : \mathbb{N} \rightarrow \mathbb{C}$ with $|a(n)|\ge 1$ for all $n \geq 1,$ that is
$$
\int_{0}^{1} \Big | \sum_{n \leq N} a(n) e(n \alpha) \Big | d \alpha \geq c \log N
$$
for an absolute constant $c > 0$.
This lower bound is remarkably general, however not much is known about the extremizers. For coefficients of modulus one, generically we expect the left-hand side to be closer to the trivial upper bound $N^{1/2}$, which follows from the Cauchy--Schwarz inequality and observation that $\int_{0}^{1} \Big | \sum_{n \leq N} a(n) e(n \alpha) \Big |^2 d \alpha=\sum_{n\le N}|a(n)|^2=N.$ In general, the upper bound is achieved in many cases, in particular when $S$ is a lacunary set, or a random positive density subset of an interval. Interestingly, the only examples we know that achieve the lower bound $O(\log N)$ are the unions of a bounded number of arithmetic progressions and a small unstructured set.

The inverse Littlewood problem was posed formally by Green in \cite{GreenICM} (see also the paper of Petridis \cite[Question 5.1]{Petridis} for a more precise formulation of the conjecture) and asks for a characterization of extremal sets $S$ for which this norm is $O(\log |S|).$  
Roughly speaking, we expect that
\begin{equation} \label{eq:l1norm}
\int_{0}^{1} \Big | \sum_{n \leq N} a(n) e(n \alpha) \Big | d \alpha
\end{equation}
is small when $a(n)$ possesses some strong ``additive'' structure. See \cite{BloomGreen}, \cite{Bedert}, \cite{Petridis}, \cite{Hanson}, \cite{GreenSanders}, \cite{Pichorides} for some of the results in this direction.

In this paper, we confine ourselves to the classical framework in analytic number theory, requiring the sequence $a(n)$ to be multiplicative.
It is generally believed that multiplicative functions such as the M\"obius function or the Liouville function possess no strong additive structure, and consequently, the $L^1$ norm \eqref{eq:l1norm} should be of the order $\gg N^{1/2 - o(1)}.$ 
For $a\equiv\mu$, lower bounds for \eqref{eq:l1norm} were obtained by Balog and Perelli \cite{BalogPerelli} and by Balog and Ruzsa \cite{BalogRuzsa1}; see also \cite{BalogRuzsa2} for the closely related case of $r$-free integers. More recently, Pandey and Radziwi\l\l{} \cite{PandeyRadziwillI} showed that for both $a\equiv\mu$ and $a\equiv\lambda$ the quantity \eqref{eq:l1norm} is $\gg N^{1/4 - o(1)},$ making connection to the low lying zeros of Riemann $\zeta$ function. For $a\equiv\lambda,$ this was improved to $\gg N^{1/4}$ by Li, Sun and Wang \cite{LSW} and to $\Omega(N^{1/3-o(1)})$ by Krause, Li, Sun and Wang \cite{KLSW}. Sharper quantitative bounds are known for other arithmetic weights, such as the indicator function of square-free numbers \cite{BrudernEtAl} and the Fourier coefficients of modular forms \cite{PandeyModular}.

Our point of departure is a recent breakthrough result by Pandey and Radziwi\l\l{} \cite{PandeyRadziwillII} (see Theorem \ref{thm:PR} below) showing that a multiplicative function whose exponential sum has a small $L^1$ norm must ``pretend" to be a Dirichlet character at large primes. In particular, they established the following result, which concerns real valued functions $f.$
\begin{theorem}[Pandey--Radziwi\l\l{} {\cite{PandeyRadziwillII}}]\label{thm:PR_0}\label{PR_0}
Let $f:\N\to[-1,1]$ be multiplicative and suppose
that
\[
 \liminf_{x\to\infty}\frac1x\sum_{n\le x}f(n)^2>0.
\]
If
\begin{equation*}
 \mathcal{N}_f(x):=\int_{0}^1\left|\sum_{n\le x}f(n)e(n\alpha)\right|d\alpha=x^{o_{x\to\infty}(1)}\qquad(x\to\infty),
\end{equation*}
then there exists a real primitive Dirichlet character $\chi$ 
such that
\begin{equation*}
 \sum_{p\le x}\frac{1-f(p)\chi(p)}{p}=o\left(\sum_{p\le x}\frac{1}{p}\right)=o(\log\log x).
\end{equation*}
\end{theorem}
Clearly, for $f:\N\to\{-1,0,1\}$ such a mean-square hypothesis is equivalent to $\sum_{p:f(p)=0}1/p<\infty$.

\subsection{Main results}

Our aim in this paper is to obtain precise structural conclusions for a fixed completely
multiplicative function from an $L^1$ bound that holds at every sufficiently
large scale. The first result in this direction characterizes the functions for
which the norm has the smallest possible logarithmic growth.

\begin{theorem}\label{thm:main}
    Let $f : \N\to\U\cup \{0\}$ be a completely multiplicative function satisfying $\sum_{p : f(p) = 0} \frac{1}{p} < \infty$ (so that $\sum_{n\le x} |f(n)|^2 \gg x$) and 
    \begin{equation}\label{eq:L1hyp}
         \mathcal N_f(x) \ll \log x
    \end{equation}
    for all $x\ge 3.$ Then $f$ is of the form $f(n) = \chi(n) n^{it},$ where $\chi$ is a Dirichlet character and $t\in\R.$
\end{theorem}
One may view the conclusion of Theorem \ref{thm:main} as a strong
manifestation of the asymptotic inverse Littlewood conjecture for the class
of multiplicative functions. Interestingly, the converse also holds, in fact, it is easy to show that a function of the form
$f(n)=\chi(n)n^{it}$ satisfies \eqref{eq:L1hyp}.
Our second main result provides an optimal improvement of Theorem \ref{thm:PR_0}, which also works for complex multiplicative functions.
\begin{theorem}\label{thm:subpolynomial}
For every $\delta\in(0,1]$ there exists $\eps_0=\eps_0(\delta)>0$, depending only on $\delta$, such that the following holds. Let $f:\N\to\C$ be a completely multiplicative function with $|f|\le1$, and suppose that
\begin{equation}\label{eq:l2-density}
 \sum_{n\le x}|f(n)|^2\ge\delta x
\end{equation}
for all sufficiently large $x$. If
\begin{equation}\label{eq:subpoly-hyp}
 \limsup_{x\to\infty}\frac{\log\mathcal N_f(x)}{\log x}\le\eps_0(\delta),
\end{equation}
then there exist a primitive Dirichlet character $\chi$ and $t\in\R$
such that
\begin{equation}\label{eq:subpoly-conclusion}
 \sum_p\frac{|f(p)-\chi(p)p^{it}|}{p}<\infty.
\end{equation}
\end{theorem}
The conclusion of Theorem \ref{thm:subpolynomial} is significantly stronger than the statement that $f$ is merely ``pretentious" (that is $\mathbb{D}(f,n^{it}\chi;\infty)<\infty$). Moreover, the hypothesis \eqref{eq:subpoly-hyp} holds, in particular, if $\mathcal N_f(x)\ll_f x^{\eps_0(\delta)}$ for all $x\ge3$ and crucially the implied constant may depend on $f$ in an arbitrary way since it only affects the size of the sum in \eqref{eq:subpoly-conclusion} and not the constant $\eps_0(\delta)$. Finally, if $f:\N\to\U\cup\{0\}$ is completely multiplicative and $\sum_{p:f(p)=0}1/p<\infty$, then clearly
\[
 \lim_{x\to\infty}\frac1x\sum_{n\le x}|f(n)|^2=\prod_{p:f(p)=0}\Bigl(1-\frac1p\Bigr)>0,
\]
so Theorem~\ref{thm:subpolynomial} applies with $\delta=\frac12\prod_{p:f(p)=0}(1-1/p),$ say. 

It is natural to ask if the hypothesis of Theorem \ref{thm:subpolynomial} are necessary, in particular, whether one can seek similar characterization if $\limsup$ in \eqref{eq:subpoly-hyp} is replaced by a weaker information on various scales.
In this direction, Theorem~\ref{thm:nonpretentious-example} below shows that $\limsup$ in \eqref{eq:subpoly-hyp} cannot be replaced by $\liminf,$ and so information about the behavior of $\mathcal{N}_f(x)$ on {\it all} large scales $x\ge 3$ is required to obtain meaningful rigidity results.
\begin{theorem}\label{thm:nonpretentious-example}
There is a completely multiplicative function $f:\N\to\{-1,1\}$
which is non-pretentious and satisfies
\begin{equation}\label{eq:example-uniform}
 \mathcal N_f(x)=o(\sqrt{x})\qquad(x\to\infty).
\end{equation}
Moreover,
\begin{equation}\label{eq:example-exponents}
 \liminf_{x\to\infty}\frac{\log\mathcal N_f(x)}{\log x}=0,
 \qquad
 \limsup_{x\to\infty}\frac{\log\mathcal N_f(x)}{\log x}=\frac12.
\end{equation}
\end{theorem}
Curiously, this example also supplies cancellation $o(\sqrt{x})$, which is barely below the typical size of $O(\sqrt{x}),$ which is attained, for example for a random (Rademacher or Steinhaus) multiplicative functions. An interested reader may surely derive this from the work of Benatar, Nishry and Rodgers \cite{BNR} or more directly from S. Hardy \cite{H}.
Unfortunately, we were unable to determine an optimal value $\varepsilon_0(\delta)=\varepsilon_0(1)$ for completely multiplicative functions $f:\mathbb{N}\to \{-1,1\}.$ It seems plausible that 
\[ \limsup_{x\to\infty}\frac{\log\mathcal N_f(x)}{\log x}=\frac12,\]
which remains an interesting open problem. One may regard this problem as an $L_1$ version of a well-known conjecture of Sark\"{o}zy, claiming that for every completely multiplicative $f:\mathbb{N}\to \{-1,1\},$
\[\limsup_{x\to\infty}\frac{\|S_{f}(x,\alpha)\|}{\sqrt{x}}=\infty.\]
\subsection{Notations.}\label{pretentious}

For $1$-bounded multiplicative functions $f,g$ and for $1\le y<x$, define the pretentious distance
\[
 \D(f,g;y,x)^2
 :=\sum_{y<p\le x}\frac{1-\re(f(p)\overline{g(p)})}{p}.
\]
We also write
\[
 S_f(x,\alpha):=\sum_{n\le x}f(n)e(\alpha n).
\]
The iterated logarithms are defined by $\log_1x=\log x$ and $\log_{k+1}x=\log(\log_kx)$, whenever the right-hand side is defined. We use $\logstar x$ for the least $k$ such that $\log_k x\le 10$.

\subsection{Outline of the proofs}
One of the goals of the present work is to describe a general set of ideas, what we call "pretentious boosting," which, roughly speaking, allows to upgrade the ``moderately pretentious" bound  $\mathbb{D}(f,\chi(n)n^{it};X)=o(\log\log X)$ (see \cite{KMT} for the terminology) to a much stronger conclusion $\mathbb{D}(f,\chi(n)n^{it};\infty)=O(1)$ in various problems. 
Some of these techniques have already been used in the works of the first author \cite{Bazin}, of the first author with Pylaiev and Tyrrell \cite{BPT} and of the second author with Mangerel and Ter\"av\"ainen \cite{KMT}. It is worth mentioning that the proofs of Theorems~\ref{thm:main} and Theorem~\ref{thm:subpolynomial} follow different strategies, and while it is possible to unify the approaches, we opted not to do so in order to illustrate different and complementary ideas (in fact, we shall use them in the future works). 

Our common point of departure is the following important theorem of Pandey and
Radziwi\l\l~\cite{PandeyRadziwillII}, which associates a character and an
Archimedean twist to a multiplicative function whose exponential sum has
small $L^1$ norm.

\begin{theorem}[Pandey--Radziwi\l\l{} {\cite[Main Theorem~A]{PandeyRadziwillII}}]\label{thm:PR}\label{PR}
Let $f$ be a $1$-bounded multiplicative function, let $c\in(0,1)$, and let $x,\Delta\ge1/c$ be such that
\[
 \sum_{n\le x}|f(n)|^2\ge cx
 \qquad\text{and}\qquad
 \int_0^1 |S_f(x,\alpha)|\,d\alpha\le \Delta.
\]
Then there are a primitive Dirichlet character $\chi$ of conductor $m$ and $t\in\R$, with
\[
 m(1+|t|)\ll_c \Delta^2,
\]
such that
\[
 \D(f,\chi(n)n^{it};\Delta^2,x)\ll_c1.
\]
\end{theorem}
Note that, in \cite{PandeyRadziwillII} the dependence of the implied constants on $c$ is made rather explicit, and the sum over primes defining the distance is restricted to $(p,m)=1,$ however this restriction changes the squared distance only by $O_c(\Delta^{-2})$, since $m\ll_c\Delta^2$. In what follows, we shall only use that the implied constants depend on $c$ alone.

The first step is then to show that the character and
Archimedean twist supplied by Theorem \ref{thm:PR} are independent of $x,$ provided we have uniform control of $L^1$ norm for all $x\ge 3.$
This is, by now, standard stabilization argument, which is executed in Corollaries~\ref{cor:stabilization} and~\ref{cor:stabilization-fixed}.\\
We now describe a general underlying principle of our proofs:
if distance is unusually small on many scales, say $\D(f,\chi(n)n^{it};\Delta^2,x)\ll 1,$ then by pigeonholing one may locate many scales $[y(X),X]$ on which  $\D(f,\chi(n)n^{it};y(X),X)=o_{X\to\infty}(1),$ for appropriately chosen growing $y(X)\to\infty.$ The key observation now is that on such scales $[X,2X],$ the contribution of large primes $p>Y(x)$ can be well controlled, in particular using ``pretentious" methods developed over the past decade.
A similar idea was also used in the work of the second author with Mangerel and Ter\"av\"ainen \cite{KMT}, where this principle allowed to establish  Elliott conjecture along dense subsecuence of scales.

In order to prove Theorem \ref{thm:main}, we adapt a uniform Delange-type
asymptotic yielding large values of exponential sum at carefully chosen rational frequencies, with explicit Euler
factors recording the discrepancies between $f(p)$ and $\chi(p)p^{it}.$ This is a content of Proposition \ref{prop:majorarc} where all the error terms crucially depend on $\mathbb{D}(f,n^{it}\chi;y(X),X)\to 0$ when $X\to\infty.$
It is then necessary to produce sufficiently large ``major arcs" to lower bound the $L^1-$norm. To this end, we use simple albeit important averaging argument in Proposition \ref{prop:averaging} turning large exponential sums at rational
points into lower bounds for the $L^1$ norm on large surrounding arcs, yielding the desired contradiction provided $f(p)\ne \chi(p)p^{it}$ for some $p.$ These arguments are outlined in Sections~\ref{sec:majorarc}-\ref{sec:main-end}.\\
The key novelty in the proof of Theorem \ref{thm:subpolynomial} is the connection of $L_1$ norm and the ``pretentious distance" via the following modified Elliott-type inequality. 
\begin{proposition}\label{elliott-intro}
%we give here the statement incorporating the smooth cutoff, since it is what we are using in section 7
Let $\chi$ be a primitive Dirichlet character of conductor $m$ and
$t\in\R$. Let $f$ be a $1$-bounded function and let $X,\Delta$ satisfy
\[
 \int_0^1\left|\sum_{n\le x}f(n)e(\alpha n)\right|\,d\alpha
 \le\Delta\qquad(1\le x\le X).
\]
Let $P\ge 1000m^2$ and $\kappa>0$ such that
$3\log P\le1/\kappa\le\log X$, and let $\Phi$ be a continuous, piecewise $C^1$ function supported on $[P^3, X]$ with $|\Phi|\le1$ and $|\Phi'(u)|\le 1/u$. We then have
\begin{equation}\label{elliott-eq-intro}
 \sum_{\substack{p\le P\\p\nmid m}}
 \left|\sum_n\frac{f(n)\Phi(n)\conj\chi(n)}{n^{1+\kappa+it}}
       -p\sum_{p\mid n}\frac{f(n)\Phi(n)\conj\chi(n)}{n^{1+\kappa+it}}\right|
 \ll (1+|t|)\frac{\Delta}{\kappa \log P},
\end{equation}
where the implied constant is absolute.
\end{proposition}
This inequality is proved in Section~\ref{sec:5} for general
$1$-bounded functions, without any use of multiplicativity. Its underlying idea is again that large exponential
sums at rationals of the form $a/(mp)$ already force a significant contribution to
the $L^1$ norm. Here the sums are  weighted by $n^{-1-\kappa-it}$, allowing the
left-hand side to be expressed in terms of a convergent Euler product of the corresponding Dirichlet series when $f$ is
completely multiplicative, after applying Corollary~\ref{cor:stabilization-fixed} selecting suitable twist $\chi(n)n^{it}$.\\
In order to appreciate the role of this smooth cut-off, we remark that the same simple proof as in Proposition \ref{elliott-intro} would immediately give the bound
\begin{equation}\label{elliott-eq1}
 \sum_{\substack{p\le P\\p\nmid m}}
 \left|\sum_n\frac{f(n)\conj\chi(n)}{n^{1+\kappa+it}}
       -p\sum_{p\mid n}\frac{f(n)\conj\chi(n)}{n^{1+\kappa+it}}\right|
 \ll \frac{P\kappa\log P+PX^{-\kappa}}{\kappa\log P}+\frac{\Delta}{\kappa\log P},
\end{equation}
where crucially the first error term arises from estimating sharp cut-off trivially. We can then explicitly evaluate the quantity
\begin{equation*}
 \begin{aligned}
 &\sum_n\frac{f(n)\conj\chi(n)}{n^{1+\kappa+it}}
   -p_0\sum_{p_0\mid n}\frac{f(n)\conj\chi(n)}{n^{1+\kappa+it}}\\
 &\hspace{2em}=
 \left(1-\frac{f(p_0)\conj\chi(p_0)}{p_0^{\kappa+it}}\right)
 L(f\conj\chi,1+\kappa+it),
 \end{aligned}
\end{equation*}
and use the standard bound
\[
 |L(f\conj\chi,1+\kappa+it)|
 \gg\frac1\kappa
       \exp\{-\D^2(f,\chi(n)n^{it};X)\}.
\]
to obtain a satisfactory upper bound on $\mathbb{D}(f,1,\chi(n)n^{it},P)$ via the $L_1$ control in \eqref{elliott-eq1}. The straightforward calculation then shows that iterating such a bound at the scales $X_k=2^k$ already yields the desired conclusion of Theorem \ref{thm:subpolynomial} in a wide range $\Delta=\Delta_X=O(\exp((\log X)^{1-\varepsilon} ).$\\

In order to improve this range further, we use an explicit smoothed cutoff $\Phi$ localizing the left hand side of \eqref{elliott-eq-intro} to the intervals $[P^3,P^A]$ for sufficiently large $A.$
This allows us to cancel the trivial contributions of pole at $\kappa=0$ (which was present in both terms of \eqref{elliott-eq1}) and furthermore control the distances $\mathbb{D}(f,\chi(n)n^{it};P^A)$ at scale $P^A$ via the distance $\mathbb{D}(f,\chi(n)n^{it};P)$ (see Proposition \ref{prop:direct-elliott}). We then use the iterative argument and partial summation to get a uniform bound of Theorem \ref{thm:subpolynomial}. Note that all constants in these arguments depend on $f$ only through the constant $\delta$ in \eqref{eq:l2-density}, which is the main reason why $\eps_0$ in Theorem~\ref{thm:subpolynomial} depends only on density $\delta$.

Finally, in
Section~\ref{sec:example} we prove Theorem \ref{thm:nonpretentious-example}, where we give the nonpretentious example, which
does not satisfy any fixed power bound $x^a$ with $a<1/2$ and behaves ``wildly" on various scales. This fairly simple construction is accomplished by taking $f(p)=-1$ on lacunary but long blocks of primes.

\section{``Pretentiousness" is stable}\label{sec:pretentious}
Our first task is to select the character and the Archimedean twist
independently of the scale after applying Theorem \ref{thm:PR}. This is a fairly standard task by now and appeared in many work on the subject over the past decade (see \cite{KluComp}, \cite{MangerelMathannalen}, \cite[\S13]{PandeyRadziwillII}). The following statement provides somewhat more general control over the distance function suitable for our applications in Theorems~\ref{thm:main} and~\ref{thm:subpolynomial}.

\begin{corollary}\label{cor:stabilization}\label{PR-uniform}
Let $c > 0$ and $f$ be a $1$-bounded multiplicative function with
$\sum_{n\le x}|f(n)|^2\ge cx$ for all sufficiently large $x$. Suppose that
\[
 \int_0^1|S_f(x,\alpha)|\,d\alpha\le x^{\psi(x)}
\]
for all sufficiently large $x$, where $\psi(x)\ge0$ and $\psi(x)\to0$.
Then there are a primitive Dirichlet character $\chi$ and $t\in\R$ such
that, for all sufficiently large $x$,
\begin{equation}\label{eq:stabclaim-subpoly}
 \D(f,\chi(n)n^{it};x^{2\psi(x)},x)^2\ll 1.
\end{equation}
Moreover, if
$\int_0^1|S_f(x,\alpha)|\,d\alpha\ll\log x$, then the same $\chi,t$ satisfy
\begin{equation}\label{eq:stabclaim}
 \D(f,\chi(n)n^{it};(\log x)^2,x)^2\ll 1,
\end{equation}
and if $\psi(x)\ll(\log x)^{-\eps}$ for some $\eps>0$, then for $x\ge 2$
\begin{equation}\label{eq:stabclaim-global}
 \D(f,\chi(n)n^{it};x)^2\ll \log_3x.
\end{equation}
\end{corollary}

\begin{proof}
The proof is a straightforward adaptation of arguments from
\cite[Lemmas~12.2, 12.3 and 12.5]{PandeyRadziwillII} and so we shall only briefly comment on the required changes.
We begin with the following important ``repulsion" claim: for each fixed $C_0>0$, there are constants $A_0,C_1>0$, depending only on $C_0$, with the
following property. Suppose that $\chi_i$ are primitive characters of
conductors $q_i$, that $q_i(1+|t_i|)\le \Delta_0$ for $i=1,2$, such that
\[
 y\ge4\Delta_0^2,\qquad X\ge y^{A_0},\qquad
 \D(\chi_1(n)n^{it_1},\chi_2(n)n^{it_2};y,X)^2\le C_0.
\]
Then, for sufficiently large $X$,
\begin{equation}\label{eq:twist-separation}
 \chi_1=\chi_2,\qquad |t_1-t_2|\le C_1/\log X.
\end{equation}
To see this, put $u=t_1-t_2$ and let $\xi$ be the primitive character
inducing $\chi_1\overline{\chi_2}$. Its conductor $q$ satisfies the bound
\[
 q(1+|u|)
 \le q_1q_2(1+|t_1|+|t_2|)
 \le y,
\]
 and $\D(1,\xi(n)n^{iu};y,X)^2\le C_0$. The lemmas in \cite[\S12]{PandeyRadziwillII} are stated for real-valued functions and below they are only applied to the real-valued functions $1$ and $\xi$ (the latter once $\xi$ is known to be real), which is legitimate since $q\le y$ and $|u|\le y$. If $\xi$ is not real, Lemma~12.2 of \cite{PandeyRadziwillII}, applied with $f\equiv1$, gives $\D(1,\xi(n)n^{iu};y,X)^2\gg\log A_0$ and consequently $\xi$ is real provided $\log A_0\gg C_0$ is large. Applying Lemma~12.3 of \cite{PandeyRadziwillII} to the real character $\xi$, for $|u|\log X\ge 10,$ say we get an immediate lower bound
\[
 \D(1,\xi(n)n^{iu};y,X)^2\gg\log\min\{A_0,|u|\log X\}.
\]
Consequently, provided $\log A_0\gg C_0$, we conclude that $|u|\log X\le C_1$ for some $C_1.$
The triangle inequality and the bound
\[
 \D(1,n^{iu};X)^2=\sum_{p\le X}\frac{1-\cos(u\log p)}p
 \le\frac{u^2}2\sum_{p\le X}\frac{(\log p)^2}p\le C_1^2
\]
now give $\D(1,\xi;y,X)\le C_0^{1/2}+C_1$.
If $\xi$ were non-principal, the proof of Lemma~12.5 of \cite{PandeyRadziwillII}, applied to the
real characters $1$ and $\xi$, would give the lower bound
$\D(1,\xi;y,X)^2\gg\log A_0,$ which yields a contradiction provided $A_0$ is again chosen large enough. Thus $\xi$ is
principal and $\chi_1=\chi_2,$ which establishes \eqref{eq:twist-separation}. 

From now on, all constants are allowed to depend on $c$. We fix a constant $\Delta_0\ge1$ large enough and put
\[
 \Delta_x=\max\left\{\Delta_0,
                 \int_0^1|S_f(x,\alpha)|\,d\alpha\right\}.
\]
By our assumption $\Delta_x=x^{o(1)}$ and so by Theorem~\ref{thm:PR}, for each
sufficiently large $x$ we may choose a primitive character $\chi_x$,
of conductor $q_x$, and $t_x\in\R$ such that
\begin{equation}\label{eq:uniform-model-at-x}
 \begin{aligned}
 q_x(1+|t_x|)&\ll \Delta_x^2,\\
 \D(f,\chi_x(n)n^{it_x};\Delta_x^2,x)^2&\le C.
 \end{aligned}
\end{equation}

We now take a sequence $x_k=\exp(2^k)$ and select $A_0$ in
\eqref{eq:twist-separation} and moreover put
$y_k=x_k^{1/A_0}$. We note that
$\Delta_{x_k}^2,\Delta_{x_{k+1}}^2\le y_k$ for all
sufficiently large $k$ and by the triangle inequality on $(y_k,x_k]$, followed by
\eqref{eq:twist-separation}, we deduce
\begin{equation}\label{eq:stabilization-increments}
 \chi_{x_k}=\chi_{x_{k+1}},\qquad
 |t_{x_k}-t_{x_{k+1}}|\ll\frac1{\log x_k}.
\end{equation}
Thus the characters are eventually equal to a fixed primitive character
$\chi.$ Since $\sum_k1/\log x_k<\infty$ we have $t_{x_k}\to t,$ when $k\to\infty$ with 
\[
 |t_{x_k}-t|\ll\frac1{\log x_k}.
\]

For arbitrary $x_k\le x<x_{k+1}$, we apply the same arguments as before, noting that
$x_{k+1}\le x^2$ and $\Delta_x=x^{o(1)}$ uniformly for all sufficiently
large $x$ to arrive at
\begin{equation}\label{eq:stabilization-rate}
 \chi_x=\chi,\qquad |t_x-t|\ll\frac1{\log x}.
\end{equation}
In particular,
\[
 \begin{split}
 \D(\chi(n)n^{it_x},\chi(n)n^{it};\Delta_x^2,x)^2
 &\ll(t_x-t)^2\sum_{p\le x}\frac{(\log p)^2}p+O(1)
 \ll 1.
 \end{split}
\]
Combining this with \eqref{eq:uniform-model-at-x} and the triangle
inequality yield
\begin{equation}\label{eq:uniform-model-limit}
 \D(f,\chi(n)n^{it};\Delta_x^2,x)^2\ll1.
\end{equation}

Since $\Delta_x\le\Delta_0x^{\psi(x)}$, replacing this in
\eqref{eq:uniform-model-limit} by $x^{2\psi(x)}$ contributes at most
$O(1)$, proving \eqref{eq:stabclaim-subpoly}. Finally, if $\Delta_x\ll\log x,$ the same comparison
proves \eqref{eq:stabclaim}, for the same choice of $\chi,t$ as before.
The bound \eqref{eq:stabclaim-global} now follows by applying \eqref{eq:stabclaim-subpoly} iteratively at the scales $e^{x}$ and $e^{x^{1-\varepsilon}}$ for $0<\varepsilon<1$ after $O_\eps(\log_3x)$ iterations. 

\end{proof}

For the proof of Theorem~\ref{thm:subpolynomial} we need a variant of Corollary~\ref{cor:stabilization} in which $\mathcal N_f(x)$ is only assumed to be bounded by a fixed small power of $x$, and in which all constants depend on $f$ only through its $\ell^2$-norm. The proof is essentially the same as of the previous Corollary~\ref{PR-uniform}, so we will remain brief.

\begin{corollary}\label{cor:stabilization-fixed}
For every $\delta\in(0,1],$ there are constants $\eps_1=\eps_1(\delta)>0$ and $\mathcal D=\mathcal D(\delta)\ge1$, depending only on $\delta$, with the following property. Let $f$ be a $1$-bounded multiplicative function and suppose that, for all sufficiently large $x$,
\[
 \sum_{n\le x}|f(n)|^2\ge\delta x
 \qquad\text{and}\qquad
 \mathcal N_f(x)\le x^{\eps+\eta(x)},
\]
where $0<\eps\le\eps_1$, $\eta(x)\ge0$ and $\eta(x)\to0$ as $x\to\infty$. Then there are a primitive Dirichlet character $\chi$ and $t\in\R$ such that, for all sufficiently large $x$,
\begin{equation}\label{eq:stabclaim-fixed}
 \D(f,\chi(n)n^{it};x^{3\eps},x)^2\le\mathcal D.
\end{equation}
\end{corollary}

\begin{proof}
The proof is essentially the same as of the previous lemma. Let $C_\delta\ge1$ be such that Theorem~\ref{thm:PR}, applied with $c=\delta/2$, holds in the form $m(1+|t|)\le C_\delta\Delta^2$ and $\D(f,\chi(n)n^{it};\Delta^2,x)^2\le C_\delta$. Let $A_0$ and $C_1$ be the constants of \eqref{eq:twist-separation} corresponding to $C_0:=4C_\delta$, and put
\[
 \eps_1:=\frac1{20A_0},\qquad \mathcal D:=2C_\delta+20C_1^2.
\]
Let $\eta^*(x):=\sup_{u\ge x}\eta(u)$, so that $\eta^*(x)\to0$. For large $x$ we take $\Delta_x:=\max\{2/\delta,x^{\eps+\eta(x)}\}\ge\mathcal N_f(x)$, and let $\chi_x$, of conductor $q_x$, and $t_x\in\R$ be given by Theorem~\ref{thm:PR}, so that
\[
 q_x(1+|t_x|)\le C_\delta\Delta_x^2,
 \qquad
 \D(f,\chi_x(n)n^{it_x};\Delta_x^2,x)^2\le C_\delta.
\]
We now consider $x\le x'\le x^2$ and set $y:=x^{1/A_0}$. Both pairs $(\chi_x,t_x)$ and $(\chi_{x'},t_{x'})$ satisfy $q(1+|t|)\le\Delta_0:=C_\delta x^{4(\eps+\eta^*(x))}$ and since $8\eps A_0\le2/5$, for all sufficiently large $x$ we have $4\Delta_0^2\le y$ and $\Delta_{x'}^2\le y$, so that both distance bounds hold on $(y,x]$. The triangle inequality now gives
\[
 \D(\chi_x(n)n^{it_x},\chi_{x'}(n)n^{it_{x'}};y,x)^2\le 4C_\delta=C_0,
\]
and \eqref{eq:twist-separation} yields $\chi_{x'}=\chi_x$ and $$|t_{x'}-t_x|\le C_1/\log x.$$ Applying this again with $x_k=\exp(2^k)$ and $x'=x_{k+1}=x_k^2$, we see that the characters $\chi_{x_k}$ are eventually equal to a fixed primitive character $\chi$, of conductor $m$ say, and that $t_{x_k}\to t$ with
\[
 |t_{x_k}-t|\le\sum_{j\ge k}\frac{C_1}{\log x_j}=\frac{2C_1}{\log x_k}.
\]
For $x_k\le x\le x_{k+1}=x_k^2$ we apply the same argument to the pair $(x_k,x)$ and obtain $\chi_x=\chi$ and $|t_x-t|\le 6C_1/\log x$. Hence, for large $x$, we have
\[
 \D(\chi(n)n^{it_x},\chi(n)n^{it};\Delta_x^2,x)^2
 \le\frac{(t_x-t)^2}2\sum_{p\le x}\frac{(\log p)^2}p
 \le10C_1^2,
\]
and the triangle inequality gives $\D(f,\chi(n)n^{it};\Delta_x^2,x)^2\le(C_\delta^{1/2}+\sqrt{10}\,C_1)^2\le\mathcal D$. Finally, we observe that $\Delta_x^2= x^{2\eps+2\eta(x)}\le x^{3\eps}$ for all large $x$, since $\eta(x)\to0,$ which completes the proof of \eqref{eq:stabclaim-fixed}.
\end{proof}

\section{A Delange-type asymptotic and major arcs}\label{sec:majorarc}

We next establish a rather standard uniform Delange-type asymptotic for the mean values of multiplicative functions (Lemma~\ref{lem:delange}) and then naturally extend it to a similar result for the exponential sums (Proposition~\ref{prop:majorarc}). Crucially, the error terms depend on the behaviour of $f$ on large primes only. 
We note that Lemma~\ref{lem:delange} could have easily been deduced from a more general correlation statement in \cite{KMT} or following the proof in \cite{KluComp} with somewhat stronger error terms. Here we instead deduce a cheaper version, which follows from the quantitative Kubilius model in \cite[Th\'eor\`eme~1, (1.7) and (1.13)]{TenenbaumKubilius}. 
\begin{lemma}\label{lem:delange}
Let $t\in\R$ be fixed and let $f:\N\to\C$ be completely multiplicative $|f|\le 1.$ Let
 $F(n):=f(n)n^{-it},$
  and $y:=x^{1/\log_2x}$ ,
and put
\[
 \mathbb{D}(F,1;y,x)^2=D_y^2:=\sum_{y<p\le x}\frac{1-\re F(p)}p.
\]
Then
\begin{equation}
 \sum_{n\le x}f(n)
 =\frac{x^{1+it}}{1+it}
 \prod_{p\le x}\frac{p-1}{p-f(p)p^{-it}}
 +O_t\left(xD_y\sqrt{\log_3x}+x\exp(-c\log_2x\log_3x)\right),\label{eq:delange}
\end{equation}
where $c>0$ is absolute.
\end{lemma}

\begin{proof}
Define the completely multiplicative function $F_y$ by
\[
 F_y(p)=
 \begin{cases}
 F(p),&p\le y,\\
 1,&p>y,
 \end{cases}
 \qquad
 P_y:=\prod_{p\le y}\frac{p-1}{p-F(p)}.
\]
We use the quantitative Kubilius model in
\cite[Th\'eor\`eme~1, (1.7) and (1.13)]{TenenbaumKubilius},
setting parameter $\varepsilon=1/2$. We have an asymptotic formula
\[
 \sum_{n\le z}F_y(n)
 =zP_y+O\left(z u_z^{-u_z}+z^{1/2}\right),
 \qquad u_z=\frac{\log z}{\log y},\quad z\ge y\ge2,
\]
 For
$x^{1/2}\le z\le x$, with these choices of parameters we have for some absolute constant $c>0,$
\[
 \sum_{n\le z}F_y(n)
 =zP_y+O\left(z\exp(-c\log_2x\log_3x)\right).
\]
Applying partial summation gives
\[
 \sum_{n\le x}F_y(n)n^{it}
 =\frac{x^{1+it}}{1+it}P_y
  +O_t\left(x\exp(-c\log_2x\log_3x)\right),
\]
where the sum of the terms below $x^{1/2}$ trivially contributes $O(x^{1/2}).$
For every $n$, telescoping the factors at primes $p>y$ gives
\[
 |F(n)-F_y(n)|
 \le\sum_{p>y}v_p(n)|F(p)-1|,
\]
where $v_p(n)$ is a $p-$valuation of $n$. Therefore
\[
 \sum_{n\le x}|F(n)-F_y(n)|
 \le x\sum_{y<p\le x}\frac{|F(p)-1|}{p-1}.
\]
By Cauchy's inequality, $|F(p)-1|^2\le2(1-\re F(p))$, and
Mertens' theorem,
\[
 \sum_{y<p\le x}\frac{|F(p)-1|}{p-1}
 \ll D_y\left(\sum_{y<p\le x}\frac1p\right)^{1/2}
 \ll D_y\sqrt{\log_3x}.
\]
Thus we can replace $F_y(n)n^{it}$ by $f(n)$ at acceptable cost
$O(xD_y\sqrt{\log_3x})$.

Finally, each factor $(1-1/p)/(1-F(p)/p)$ has modulus at most $1$,
and
\[
 \left|\frac{1-1/p}{1-F(p)/p}-1\right|
 =\frac{|F(p)-1|}{|p-F(p)|}
 \le\frac{|F(p)-1|}{p-1}.
\]
Another telescoping argument therefore gives
\[
 \left|
 \prod_{p\le x}\frac{1-1/p}{1-F(p)/p}-P_y
 \right|
 \le\sum_{y<p\le x}\frac{|F(p)-1|}{p-1}
 \ll D_y\sqrt{\log_3x}.
\]
Combining these estimates establishes \eqref{eq:delange}.
\end{proof}

We now apply the preceding lemma to obtain an asymptotic formula at the rational multiplies of $m.$
To this end, we define
\begin{equation}
 P_\chi(f,x):=\prod_{p\le x}\frac{p-1}{p-f(p)\overline{\chi(p)}p^{-it}}.\label{eq:Pchi}
\end{equation}

\begin{proposition}\label{prop:majorarc}
Let $f:\N\to\C$ be completely multiplicative and $1$-bounded. Let $\chi$ be a fixed primitive Dirichlet character of conductor $m$, and let $t\in\R$ be fixed. Suppose that
\begin{equation}
 \D(f,\chi(n)n^{it};x^{1/\log_2x},x)^2
 \ll\frac{\log_3x}{\log_2x}\label{eq:ma-hyp1}
\end{equation}
and
\begin{equation}
 \D(f,\chi(n)n^{it};1,x)^2\ll\log_3x.\label{eq:ma-hyp2}
\end{equation}
Let $q\le\log x$ and let $(a,qm)=1$. Then
\begin{equation}
 S_f\left(x,\frac{a}{qm}\right)
 =\frac{x^{1+it}}{1+it}
 \frac{\overline{\chi(a)}\tau(\chi)}{\varphi(m)}
 P_\chi(f,x)\mathcal E_q(f,\chi,t)
 +O_{m,t,\eps}\left(\frac{x}{(\log_2x)^{1/2-\eps}}\right),\label{eq:majorarcformula}
\end{equation}
where the Gauss sum
\[
 \tau(\chi):=\sum_{b\bmod m}\chi(b)e(b/m),
\]
and $\mathcal E_q$ is multiplicative in $q$, with local factors
\begin{equation}
 \mathcal E_{p^\nu}(f,\chi,t)
 =
 \begin{cases}
 \displaystyle
 \frac{f(p)^{\nu-1}\bigl(f(p)-\chi(p)p^{it}\bigr)}{\varphi(p^\nu)p^{\nu it}},& p\nmid m,\\[2.2ex]
 \displaystyle
 \frac{f(p)^\nu}{p^{\nu(1+it)}},& p\mid m.
 \end{cases}\label{eq:localfactor}
\end{equation}
\end{proposition}

\begin{proof}
We write $Q=qm$ and for a character $\psi\bmod k$, let
\[
 G_k(\psi):=\sum_{\substack{b\bmod k\\ (b,k)=1}}e(b/k)\overline{\psi(b)}.
\]
Now if $(r,k)=1$, we have
\begin{equation}
 e(ar/k)=\frac1{\varphi(k)}\sum_{\psi\bmod k}\psi(a)G_k(\psi)\psi(r),\label{eq:addchar}
\end{equation}
and by complete multiplicativity, we can write
\begin{equation}
 S_f\left(x,\frac aQ\right)
 =\sum_{dk=Q}\frac{f(d)}{\varphi(k)}
 \sum_{\psi\bmod k}\psi(a)G_k(\psi)
 \sum_{n\le x/d}f(n)\psi(n).\label{eq:Sdecomp}
\end{equation}

We begin with those $\psi\bmod k$ not induced by $\overline\chi$.
Then $\xi=\chi\psi$ is non-principal and has a modulus
$r=\operatorname{lcm}(m,k)\ll\log x$. For $z=x/d$ and
$T=(\log x)^4$ we take $\sigma=1+1/\log z$ and $|u|\ll T$, the standard estimates yield
\[
 |L(\sigma+iu,\xi)|\ll\log(r(1+|u|))\ll\log_2x.
\]
 Hence, uniformly for $|u|\le T$,
\[
 \D(\chi(n)n^{it},\overline\psi(n)n^{iu};z)^2
 \ge\log_2x-O(\log_3x).
\]
Together with \eqref{eq:ma-hyp2} and the triangle inequality, this implies
\[
 M:=\min_{|u|\le T}\D(f\psi,n^{iu};z)^2
 \ge\log_2x-O((\log_2x\log_3x)^{1/2}).
\]
We now apply the quantitative form
of Hal\'asz's theorem \cite[Theorem~2.1]{BGS} to end up with
\begin{equation}
 \left|\sum_{n\le z}f(n)\psi(n)\right|
 \ll z\bigl((1+M)e^{-M}+T^{-1/2}\bigr)
 \ll \frac{z}{\log x}
       \exp\left(O((\log_2x\log_3x)^{1/2})\right).\label{eq:halasznonmain}
\end{equation}
Using $|G_k(\psi)|\le\sqrt k$ and $Q\ll \log x$ the total contribution of these
characters to \eqref{eq:Sdecomp} is
\[
 \ll\frac{x}{\log x}
 \exp\left(O((\log_2x\log_3x)^{1/2})\right)
 \sum_{dk=Q}\frac{\sqrt k}{d}
 \ll\frac{x}{(\log x)^{1/2-o(1)}}.
\]

We now turn to the characters induced by $\overline\chi$. If
$\psi\bmod k$ is induced by $\overline\chi$, then $k=m\ell$ and
the relation $dk=qm$ becomes $d\ell=q$. We write
$\psi_\ell\bmod m\ell$ for the character induced by $\overline\chi$.
Recall that
\begin{equation}
 G_{m\ell}(\psi_\ell)=\mu(\ell)\chi(\ell)\tau(\chi),\label{eq:imprimitiveGauss}
\end{equation}
with the convention that the right-hand side is zero when $(\ell,m)>1$ and $\psi_\ell(a)=\overline{\chi(a)}$.

For $p\nmid\ell$ we have $\psi_\ell(p)=\overline{\chi(p)}$,
whereas for $p\mid\ell$ we have $\psi_\ell(p)=0$. Since
$\ell\le q\le\log x$ and $x^{1/\log_2x}\gg\log x$,
Lemma~\ref{lem:delange}, applied to $f\psi_\ell$, gives the asymptotic
\begin{equation}
 \sum_{n\le x/d}f(n)\psi_\ell(n)
 =\frac{(x/d)^{1+it}}{1+it}
 P_\chi(f,x)
 \prod_{p\mid\ell}\frac{p-f(p)\overline{\chi(p)}p^{-it}}{p}
 +O_{m,t,\eps}\left(\frac{x/d}{(\log_2x)^{1/2-\eps}}\right)\label{eq:maincharsum}
\end{equation}
uniformly for $d\ell=q,$ since
replacing the Euler-product cutoff $z=\frac{x}{d}$ for $d\ll\log x$ by $x$  introduces
only an error $O(\log_2x/\log x)$.
Moreover, the errors in \eqref{eq:maincharsum} give acceptable contribution since the sum
\[
 \sum_{d\ell=q}\frac1{d\varphi(\ell)}
 =\prod_{p^a\Vert q}\frac{1+a p/(p-1)}{p^a}\le\frac32.
\]
Substituting \eqref{eq:imprimitiveGauss} and \eqref{eq:maincharsum}
into \eqref{eq:Sdecomp}, we obtain the main term
\[
 \frac{x^{1+it}}{1+it}\overline{\chi(a)}\tau(\chi)P_\chi(f,x)
 \sum_{d\ell=q}
 \frac{f(d)}{d^{1+it}}
 \frac{\mu(\ell)\chi(\ell)}{\varphi(m\ell)}
 \prod_{p\mid\ell}\frac{p-f(p)\overline{\chi(p)}p^{-it}}p.
\]
Since $\chi(\ell)=0$ unless $(\ell,m)=1$, this equals
\begin{equation}
 \frac{x^{1+it}}{1+it}\frac{\overline{\chi(a)}\tau(\chi)}{\varphi(m)}P_\chi(f,x)
 \sum_{d\ell=q}
 \frac{f(d)}{d^{1+it}}
 \frac{\mu(\ell)\chi(\ell)}{\varphi(\ell)}
 \prod_{p\mid\ell}\frac{p-f(p)\overline{\chi(p)}p^{-it}}p.\label{eq:qsum}
\end{equation}
Recall, that the sum in \eqref{eq:qsum} is multiplicative in $q$ and so for
$p^\nu\Vert q$ with $p\nmid m$, the local factor is
\[
 \frac{f(p)^\nu}{p^{\nu(1+it)}}
 -\frac{\chi(p)}{p-1}
 \frac{f(p)^{\nu-1}}{p^{(\nu-1)(1+it)}}
 \frac{p-f(p)\overline{\chi(p)}p^{-it}}p.
\]

For $p\mid m$, we have $\chi(p)=0$, so only the term with
$p\nmid\ell$ remains. The local factor is then
\[
 \frac{f(p)^\nu}{p^{\nu(1+it)}}.
\]
Simplifying these local factors gives \eqref{eq:localfactor} and
completes the proof of \eqref{eq:majorarcformula}.
\end{proof}
We next turn a large value at a rational point into a lower bound 
over the surrounding arc. The following simple averaging trick is
inspired by an earlier work of the first author \cite{Bazin}.
\begin{proposition}\label{prop:averaging}
Let $f:\N\to\C$ be $1$-bounded and let $\theta\in\R$. Suppose that for all $X\le x\le2X$,
\begin{equation}
 |S_f(x+X,\theta)-S_f(x,\theta)|\ge\kappa(X).\label{eq:averaginghyp}
\end{equation}
Then for every $\delta\ge X^{-1/2}$,
\begin{equation}
 \frac1X\int_X^{3X}\int_{\theta-\delta}^{\theta+\delta}|S_f(u,\alpha)|\,d\alpha\,du
 \gg \frac{\log X}{X}\kappa(X).\label{eq:averagingconclusion}
\end{equation}
\end{proposition}

\begin{proof}
Replacing $f(n)$ by $f(n)e(\theta n)$, we may assume that
$\theta=0$. Partial summation gives, for $|\alpha|\le\delta$,
\[
 \sum_{x<n\le x+X}f(n)
 =e(-\alpha(x+X))(S_f(x+X,\alpha)-S_f(x,\alpha))
\]
\[
 \qquad +2\pi i\alpha\int_x^{x+X}e(-\alpha u)(S_f(u,\alpha)-S_f(x,\alpha))\,du.
\]
Thus
\[
 \kappa(X)\le |S_f(x+X,\alpha)|+|S_f(x,\alpha)|
 +2\pi|\alpha|\int_x^{x+X}(|S_f(u,\alpha)|+|S_f(x,\alpha)|)\,du.
\]
Averaging over $x\in[X,2X]$, we obtain
\[
 X\kappa(X)\ll (1+|\alpha|X)\int_X^{3X}|S_f(u,\alpha)|\,du.
\]
Integrating over $\alpha\in[-\delta,\delta]$ gives
\[
 \frac1X\int_X^{3X}\int_{-\delta}^{\delta}|S_f(u,\alpha)|\,d\alpha\,du
 \gg \kappa(X)\int_{-\delta}^{\delta}\frac{d\alpha}{1+|\alpha|X}.
\]
The assumption $\delta\ge X^{-1/2}$ makes the last integral
$\gg(\log X)/X$, proving the proposition.
\end{proof}

\section{Completion of the proof of Theorem \ref{thm:main}}\label{sec:main-end}

We are now ready to prove Theorem \ref{thm:main} and in the following proof, all constants are allowed to depend on $f.$

\begin{proof}[Proof of Theorem \ref{thm:main}]
Firstly, we assume \eqref{eq:L1hyp}. Since
\[
 \sum_{p:f(p)=0}\frac1p<\infty,
\]
the integers divisible by no prime $p$ with $f(p)=0$ have positive lower density. Hence
\[
 \sum_{n\le x}|f(n)|^2\gg x.
\]
Corollary \ref{cor:stabilization} gives a primitive Dirichlet character $\chi$, of conductor $m$, and $t\in\R$ such that
\begin{equation}
 \D(f,\chi(n)n^{it};(\log x)^2,x)^2\ll1\label{eq:globalblocks}
\end{equation}
for all large $x$.

\subsection*{Step 1.}
We first obtain a global bound for the distance. We define recursively $x_0=x$ and
$x_{j+1}=(\log x_j)^2$ and apply \eqref{eq:globalblocks} on the
disjoint intervals $(x_{j+1},x_j].$ As a consequence this gives, for every $m,$
\begin{equation}
 \D(f,\chi(n)n^{it};1,x)^2\ll m+\log_2x_m,\label{eq:Kdistance}
\end{equation}
since the contribution of the initial primes is bounded by
$2\sum_{p\le x_m}1/p\ll \log_2x_m$. Iterating this procedure, we deduce the bound 
\begin{equation}
 \D(f,\chi(n)n^{it};1,x)^2\ll\logstar x.\label{eq:logstardistance}
\end{equation}
Furthermore,
\begin{equation}
 \log |P_\chi(f,x)|
 =-\D(f,\chi(n)n^{it};1,x)^2+O(1),\label{eq:Psize}
\end{equation}
and consequently 
\begin{equation}
 |P_\chi(f,x)|\gg \exp(-C\logstar x).\label{eq:Plowerlogstar}
\end{equation}
In particular,
\begin{equation}
 \log_3x\,\exp(-C\logstar x)\to\infty.\label{eq:log3beatslogstar}
\end{equation}

\subsection*{Step 2.} We next choose a scale at which the distance over large primes is very small and tending to zero.
From \eqref{eq:globalblocks}, we have
\[
 \D(f,\chi(n)n^{it};\exp((\log X)^{1/2}),X)^2\ll1.
\]
Set $y_0=X$ and $y_{j+1}=y_j^{1/\log_2y_j},$ with $y_j\gg\exp((\log X)^{1/2})$ and note that there are clearly at least
$\gg\log_2X/\log_3X$ such disjoint intervals. Since the distance is bounded,
pigeonholing gives arbitrarily large $Y$ such that
\begin{equation}
 \D(f,\chi(n)n^{it};Y^{1/\log_2Y},Y)^2
 \ll\frac{\log_3Y}{\log_2Y}.\label{eq:goodscale}
\end{equation}
Moreover, for every $y\in[Y,3Y]$, the latter implies
\begin{equation}
 \D(f,\chi(n)n^{it};y^{1/\log_2y},y)^2
 \ll\frac{\log_3y}{\log_2y},\label{eq:goodscaleuniform}
\end{equation}
since the extra primes in $[Y,y]$ contribute at most $O(1/\log Y)$. Together with
\eqref{eq:logstardistance}, this verifies the hypotheses of
Proposition~\ref{prop:majorarc}, uniformly for $y\in[Y,3Y]$.

\subsection*{Step 3.}
Suppose for a moment that there is a prime $p$ such that
\begin{equation}
 f(p)\ne0,
 \qquad
 f(p)\ne\chi(p)p^{it}.\label{eq:badp}
\end{equation}
We put
\[
 M_Y:=\left\lfloor \frac{\log_3Y}{3\log p}\right\rfloor,
\]
and take $\eps=1/10$ in Proposition~\ref{prop:majorarc}, say, and furthermore, for
$1\le k\le M_Y$, we take $q=p^k$. For $(a,p^km)=1$, we introduce the frequencies 
\[
 \theta_{a,k}:=\frac{a}{p^km},
\]
for which, by Proposition \ref{prop:majorarc}, uniformly for $Y\le y\le3Y$,
\begin{equation}
 S_f(y,\theta_{a,k})
 =C_{a,k,p}\left(\frac{y}{p^k}\right)^{1+it}P_\chi(f,y)
 +O_\eps\left(\frac{Y}{(\log_2Y)^{1/2-\eps}}\right),\label{eq:badpmajorarc}
\end{equation}
where $C_{a,k,p}$ is nonzero and
$|C_{a,k,p}|\asymp 1$. Indeed, for $p\nmid m$
the local factor is
\[
 \frac{f(p)^{k-1}(f(p)-\chi(p)p^{it})}{\varphi(p^k)p^{kit}},
\]
and if $p\mid m$ it is $f(p)^k/p^{k(1+it)},$ which confirms the claim above.
Since
\[
 P_\chi(f,y)=P_\chi(f,Y)(1+O(1/\log Y))
\]
uniformly for $Y\le y\le3Y$, and since the error in \eqref{eq:badpmajorarc} is negligible compared with
\[
 \frac{Y}{p^k}|P_\chi(f,Y)|
\]
by \eqref{eq:Plowerlogstar} (recall that $p^k\le(\log_2Y)^{1/3}$ and $\eps=1/10$), we may subtract the two asymptotics.
By the triangle inequality
\[
 |(y+Y)^{1+it}-y^{1+it}|
 \ge\bigl||(y+Y)^{1+it}|-|y^{1+it}|\bigr|=Y,
\]
and therefore for all $Y\le y\le2Y$,
\begin{equation}
 |S_f(y+Y,\theta_{a,k})-S_f(y,\theta_{a,k})|
 \gg\frac{Y}{p^k}|P_\chi(f,Y)|.\label{eq:badpkappa}
\end{equation}
We now introduce the major arcs via
\[
 \mathfrak M:=\bigcup_{1\le k\le M_Y}
 \bigcup_{\substack{1\le a\le p^km\\ (a,p^km)=1}}
 \left[\frac{a}{p^km}-Y^{-1/2},\frac{a}{p^km}+Y^{-1/2}\right].
\]
For large $Y$, these arcs are clearly disjoint, since all denominators are $\ll(\log_2Y)^{1/3}$.

We now apply Proposition~\ref{prop:averaging} at each $\theta_{a,k}$ and sum
over $a,k$ to end up with 
\[
 \frac1Y\int_Y^{3Y}\int_{\mathfrak M}|S_f(u,\alpha)|\,d\alpha\,du
 \gg
 \sum_{1\le k\le M_Y}\varphi(p^km)\frac{\log Y}{Y}
 \frac{Y}{p^k}|P_\chi(f,Y)|.
\]
Since $\varphi(p^km)\asymp p^k$, this gives
\begin{equation}
 \frac1Y\int_Y^{3Y}\int_{\mathfrak M}|S_f(u,\alpha)|\,d\alpha\,du
 \gg \log Y\,M_Y\,|P_\chi(f,Y)|.\label{eq:badplower}
\end{equation}
By \eqref{eq:Plowerlogstar} and \eqref{eq:log3beatslogstar}, the
right-hand side exceeds any fixed multiple of $\log Y$ for sufficiently
large scales, which in turn, contradicts the bound
\[
 \frac1Y\int_Y^{3Y}\int_0^1|S_f(u,\alpha)|\,d\alpha\,du\ll\log Y.
\]
We therefore conclude that, for
every prime $p$,
\begin{equation}
 f(p)=0\qquad\text{or}\qquad f(p)=\chi(p)p^{it}.\label{eq:dichotomy}
\end{equation}
In particular, if $p\mid m$, then $\chi(p)=0$, so \eqref{eq:dichotomy} gives $f(p)=0$.

\subsection*{Step 4.}
We are left to deal with the set
\[
 \mathcal Z:=\{p:f(p)=0,\ \chi(p)\ne0\}.
\]
For $p\in\mathcal Z$, Proposition \ref{prop:majorarc} with $q=p$ gives, uniformly for $Y\le y\le3Y$ and $(a,pm)=1$,
\begin{equation}
 S_f\left(y,\frac a{pm}\right)
 =-\frac{y^{1+it}}{1+it}
 \frac{\overline{\chi(a)}\tau(\chi)}{\varphi(m)}
 P_\chi(f,y)\frac{\chi(p)}{p-1}
 +O_\eps\left(\frac{Y}{(\log_2Y)^{1/2-\eps}}\right).\label{eq:zeropmajorarc}
\end{equation}
By \eqref{eq:dichotomy} and the assumption $\sum_{p:f(p)=0}1/p<\infty,$
we have
\[
 \D(f,\chi(n)n^{it};Y)^2\ll1,
\]
and consequently
\begin{equation}
 |P_\chi(f,Y)|\gg1.\label{eq:Pboundedbelow}
\end{equation}
We introduce $H_Y:=(\log_2Y)^{1/3},$ so that for each $p\in\mathcal Z$
with $p\le H_Y$, \eqref{eq:zeropmajorarc} gives, uniformly for
$Y\le y\le2Y$,
\begin{equation}
 \left|S_f\left(y+Y,\frac a{pm}\right)-S_f\left(y,\frac a{pm}\right)\right|
 \gg \frac{Y}{p}.\label{eq:zeropkappa}
\end{equation}
We now define the major arcs via
\[
 \mathfrak M_0:=\bigcup_{\substack{p\in\mathcal Z\\ p\le H_Y}}
 \bigcup_{\substack{1\le a\le pm\\ (a,pm)=1}}
 \left[\frac a{pm}-Y^{-1/2},\frac a{pm}+Y^{-1/2}\right]
\]
and note that these arcs are again disjoint for large $Y.$ Applying
Proposition~\ref{prop:averaging} we get the lower bound
\[
 \frac1Y\int_Y^{3Y}\int_{\mathfrak M_0}|S_f(u,\alpha)|\,d\alpha\,du
 \gg
 \sum_{\substack{p\in\mathcal Z\\ p\le H_Y}}
 \varphi(pm)\frac{\log Y}{Y}\frac{Y}{p},
\]
 which further implies
\begin{equation}
 \frac1Y\int_Y^{3Y}\int_{\mathfrak M_0}|S_f(u,\alpha)|\,d\alpha\,du
 \gg \log Y\,\#\{p\in\mathcal Z:p\le H_Y\}.\label{eq:zerocountlower}
\end{equation}
The assumed $L^1$ bound implies the left-hand side is $O(\log Y)$ and consequently
\[
 \#\{p\in\mathcal Z:p\le H_Y\}\ll1,
\]
so that our set $\mathcal Z$ is finite.

It remains to incorporate the additional zeros into the final conductor. To this end, we put
\[
 M:=m\prod_{p\in\mathcal Z}p,
\]
and let $\chi_1$ be the character modulo $M$ induced by $\chi$.
Then $\chi_1(p)=0$ precisely when $p\mid M$, while
$\chi_1(p)=\chi(p)$ for $p\nmid M$. By \eqref{eq:dichotomy}
and the definition of $\mathcal Z$, we have, for every prime $p$,
\[
 f(p)=\chi_1(p)p^{it}.
\]
and complete multiplicativity extends the identity to all $n\ge1.$ This completes the proof.
\end{proof}

\section{A modified Elliott inequality}\label{sec:5}
 We now establish an $L_1$ version of Elliott
inequality in which a small $L^1$ norm controls a sum over primes. This proposition lies in the heart of the proof of Theorem~\ref{thm:subpolynomial} and, in fact, does not require any multiplicativity assumption on $f.$

\setcounter{theorem}{1}
\begin{proposition}\label{elliott}
%we give here the statement incorporating the smooth cutoff, since it is what we are using in section 7
Let $\chi$ be a primitive Dirichlet character of conductor $m$ and
$t\in\R$. Let $f$ be a $1$-bounded function and let $X,\Delta$ satisfy
\[
 \int_0^1\left|\sum_{n\le x}f(n)e(\alpha n)\right|\,d\alpha
 \le\Delta\qquad(1\le x\le X).
\]
Assume that $P\ge 100m^2$ and $\kappa>0$ such that
$3\log P\le1/\kappa\le\log X$, and let $\Phi$ be a continuous, piecewise $C^1$ function supported on $[P^3, X]$ with $|\Phi|\le1$ and $|\Phi'(u)|\le 1/u$. We then have
\begin{equation}\label{elliott-eq}
 \sum_{\substack{p\le P\\p\nmid m}}
 \left|\sum_n\frac{f(n)\Phi(n)\conj\chi(n)}{n^{1+\kappa+it}}
       -p\sum_{p\mid n}\frac{f(n)\Phi(n)\conj\chi(n)}{n^{1+\kappa+it}}\right|
 \ll (1+|t|)\frac{\Delta}{\kappa \log P},
\end{equation}
where the implied constant is absolute.
\end{proposition}

For completely multiplicative $f$ ``close" to $\chi(n)n^{it}$, the
left-hand side of \eqref{elliott-eq} can be expressed in terms of the discrepancy
$|f(p)-\chi(p)p^{it}|$ with the standard trick introducing  
the weight $n^{-1-\kappa-it}$ to allow us the use of  absolutely convergent
Euler products after applying Proposition~\ref{elliott}.

\begin{proof}[Proof of Proposition~\ref{elliott}]

We remark, that in this proof all implied constants are absolute and we shall keep track of the dependence on $t$ explicitly.

For $p\nmid m$, put $G(\chi)=\sum_{a\bmod m}\chi(a)e(a/m)$, so
$|G(\chi)|=\sqrt m$. The Chinese remainder theorem gives
\begin{align*}
 \sum_{\substack{a\bmod mp\\(a,mp)=1}}\chi(a)e(an/(mp))
 &=\sum_{\substack{u\bmod m\\(u,m)=1}}
   \sum_{\substack{v\bmod p\\(v,p)=1}}
   \chi(up+vm)e((up+vm)n/(mp))\\
 &=G(\chi)\conj\chi(n)\chi(p)(p\Ind_{p\mid n}-1).
\end{align*}
Consequently,
\begin{equation}\label{step21}
 \begin{aligned}
 &\left|\sum_n\frac{f(n)\Phi(n)\conj\chi(n)}{n^{1+\kappa+it}}
       -p\sum_{p\mid n}\frac{f(n)\Phi(n)\conj\chi(n)}{n^{1+\kappa+it}}\right|\\
 &=\frac1{|G(\chi)|}
   \left|\sum_{\substack{a\bmod mp\\(a,mp)=1}}\chi(a)
      \sum_n\frac{f(n)\Phi(n)}{n^{1+\kappa+it}}e(an/(mp))\right|\\
 &\ll\sum_{\substack{a\bmod mp\\(a,mp)=1}}
     \left|\sum_n\frac{f(n)\Phi(n)}{n^{1+\kappa+it}}e(an/(mp))\right|.
 \end{aligned}
\end{equation}
We next compare each rational value with the integral over a surrounding
arc. For $\alpha,\delta\in\R$, partial summation gives
\[
 \sum_n\frac{f(n)\Phi(n)}{n^{1+\kappa+it}}e(\alpha n)
 =-\int_1^\infty S_f(u,\alpha+\delta)\,
 d\left(\frac{e(-\delta u)}{u^{1+\kappa+it}}\Phi(u)\right).
\]
Since $\kappa\le1$, and $\Phi$ is supported on $[P^3, X]$ with $|\Phi|\le1$ and $|\Phi'(u)|\le 1/u$, the derivative of $u\mapsto e(-\delta u)u^{-1-\kappa-it}\Phi(u)$ is $\ll(1+|t|)(|\delta|+1/u)u^{-1-\kappa}$. Taking absolute values therefore yields
\begin{align*}
 \left|\sum_n\frac{f(n)\Phi(n)}{n^{1+\kappa+it}}e(\alpha n)\right|
 &\ll(1+|t|)\int_{P^3}^X|S_f(u,\alpha+\delta)|
                 \frac{|\delta|+1/u}{u^{1+\kappa}}\,du\\
 &\ll(1+|t|)(|\delta|+1/P^3)
          \int_{P^3}^X\frac{|S_f(u,\alpha+\delta)|}{u^{1+\kappa}}\,du.
\end{align*}
Dividing the last expression by $|\delta|+1/P^3$ and integrating over $[-\theta,\theta],$ we get for every $\theta>0$
\begin{equation}\label{step22}
 (1+|t|)\int_{P^3}^X\int_{\alpha-\theta}^{\alpha+\theta}
 \frac{|S_f(u,\beta)|}{u^{1+\kappa}}\,d\beta\,du
 \gg\log(1+P^3\theta)
       \left|\sum_n\frac{f(n)\Phi(n)}{n^{1+\kappa+it}}e(\alpha n)\right|.
\end{equation}
We now take $\theta=1/(4mP^2)$ and sum over the reduced fractions $a/(mp)$,
with $p\le P$ and $p\nmid m$. 
The centers of the arcs are separated by at least $1/(mP^2)$, and hence the
arcs are disjoint. Since $P\ge100m^2$, we also have $\log(1+P^3\theta)=\log(1+P/(4m))\ge\frac12\log P$, and therefore
\begin{equation}\label{step23}
 \begin{aligned}
 &\log P\sum_{\substack{p\le P\\p\nmid m}}
   \sum_{\substack{a\bmod mp\\(a,mp)=1}}
   \left|\sum_n\frac{f(n)\Phi(n)}{n^{1+\kappa+it}}e(an/(mp))\right|\\
 &\hspace{2em}\ll(1+|t|)
   \int_{P^3}^X\int_0^1\frac{|S_f(u,\beta)|}{u^{1+\kappa}}\,d\beta\,du
   \ll(1+|t|)\frac{\Delta}\kappa.
 \end{aligned}
\end{equation}
Combining with \eqref{step21} concludes the proof.
\end{proof}
We now write $g(n) = f(n)\conj{\chi(n)}n^{-it}$ and use complete multiplicativity of $f.$ For every prime $p$ we have
\[
 \begin{split}
 &\sum_n\frac{g(n)\Phi(n)}{n^{1+\kappa}}
    -p\sum_{p\mid n}\frac{g(n)\Phi(n)}{n^{1+\kappa}}
 =\sum_n\frac{g(n)\Phi(n)}{n^{1+\kappa}}
    -p^{-\kappa} g(p)\sum_{n}\frac{g(n)\Phi(pn)}{n^{1+\kappa}}\\
 &\qquad=(1-g(p))\sum_n\frac{g(n)\Phi(n)}{n^{1+\kappa}}
    +g(p)\sum_n\frac{g(n)}{n^{1+\kappa}} \big(\Phi(n) - p^{-\kappa} \Phi(pn)\big).
 \end{split}
\]
Hence, by the triangle inequality, the left-hand side of \eqref{elliott-eq} satisfies
\begin{equation}\label{multiplicative-elliott}
 \begin{aligned}
 &\sum_{\substack{p\le P\\p\nmid m}} 
    \left|\sum_n\frac{g(n)\Phi(n)}{n^{1+\kappa}}
    -p\sum_{p\mid n}\frac{g(n)\Phi(n)}{n^{1+\kappa}}\right|
 \ge\left|\sum_n\frac{g(n)\Phi(n)}{n^{1+\kappa}}\right|
    \sum_{\substack{p\le P\\p\nmid m}} |1-g(p)|\\
 &\hspace{12em}-\sum_{\substack{p\le P\\p\nmid m}} |g(p)|
    \left|\sum_n\frac{g(n)}{n^{1+\kappa}} \big(\Phi(n) - p^{-\kappa} \Phi(pn)\big)\right|.
 \end{aligned}
\end{equation}

In the next section, we pick an explicit cutoff $\Phi$ and estimate the second (error) term in \eqref{multiplicative-elliott} in order to prove Theorem~\ref{thm:subpolynomial}.

\section{Proof of Theorem \ref{thm:subpolynomial}}\label{sec:6}

Theorem \ref{thm:subpolynomial} is obtained as a consequence of the following explicit result, which further elaborates on Proposition~\ref{elliott} when $f$ is completely multiplicative using an explicit cutoff $\Phi.$
\begin{proposition}\label{prop:direct-elliott}
    There are absolute constants $B,C\ge1$ with the following property. Let $f$ be a completely multiplicative $1$-bounded function, let $\chi$ be a primitive Dirichlet character of conductor $m$, let $t\in\R$, and take $g(n)=f(n)\conj{\chi(n)}n^{-it}$. Let $\mathcal D\ge0$ and $A\ge C e^{\mathcal D}$, let $P\ge 100m^2$ be sufficiently large, and put $X=P^{AB}$. Suppose that
    \[
     \int_0^1\left|\sum_{n\le x}f(n)e(\alpha n)\right|\,d\alpha
     \le\Delta\qquad(1\le x\le X)
    \]
    and that
    \[
     \D(f,\chi(n)n^{it};P^3,P^A)^2\le\mathcal D.
    \]
    Then
    \begin{equation}\label{eq:direct-elliott-bound}
     \begin{split}
     \sum_{p\le P} |1-g(p)|\le{}&C(1+|t|)A^2\Delta\log P+\omega(m)+\frac{C e^{\mathcal D}P}{A(\log P)^2}\\
     &+\frac{C e^{\mathcal D}P}{A\log P}
        \sum_p\frac{|1-g(p)|}{p^{1+1/(A\log P)}}
                               \bigl(1-p^{-1/(3\log P)}\bigr),
     \end{split}
    \end{equation}
   where $\omega(m)$ denotes the number of prime factors of $m.$
\end{proposition}
We note that $|1-g(p)|=|f(p)-\chi(p)p^{it}|$ for $p\nmid m,$ and $\omega(m)$ arises from estimating the contribution $p\mid m$ trivially.  Before introducing the explicit cutoff $\Phi,$ we state some preliminary definitions and standard analytic lemmas. We also note that in this section, the constants are absolute, unless stated otherwise.

\subsection{Analytic preliminaries}
For the moment, let $g$ be any completely multiplicative $1$-bounded
function and put
\begin{equation}\label{eq:subpoly-notation}
 \begin{gathered}
 F(s)=\sum_{n\ge1}\frac{g(n)}{n^{1+s}}
     =\prod_p(1-g(p)p^{-1-s})^{-1},\qquad
 G(s)=sF(s),\quad \re s>0.
 \end{gathered}
\end{equation}
We note that the series converge absolutely for $\re s>0$ and moreover $F$ and
$G$ do not vanish there. 

\begin{lemma}\label{lem:vertical-euler}
Uniformly for $0<\sigma\le1$ and $v\in\R$,
\begin{equation}\label{eq:vertical-F}
 \begin{gathered}
 |F(\sigma+iv)|\ll(1+|v|/\sigma)^2|F(\sigma)|,\\
 |F(\sigma)|\gg\sigma,\qquad
 |F(\sigma/2)|\ll|F(\sigma)|.
 \end{gathered}
\end{equation}
\end{lemma}
\begin{proof}
These are rather standard estimates and can be deduced by applying for example \cite[Lemmas~1 and~2]{MVMeanValues} say, to $F(s-1)$. For the
lower bound, we can further use straightforward Euler product computations.
\end{proof}

\begin{lemma}\label{lem:horizontal-G}
    Let $x > 3$ and $0 < \kappa \le 1/\log x.$ Then 
    \begin{equation}\label{eq:horizontal-G}
        \abs{G\bra{\kappa + \frac{1}{\log x}}} \ll |G(\kappa)|  e^{\D^2(g, 1, x,  e^{1/\kappa})}.
    \end{equation}
\end{lemma}
\begin{proof}
    Noting $\sigma = \kappa + \frac{1}{\log Y},$ we have using the Euler products and the estimate $\zeta(1+s) \asymp 1/s$ for $0 < s \le 2/\log x,$
    \[
     \log\frac{|G(\sigma)|}{|G(\kappa)|}
     = O(1)+\sum_p\frac{1-\re g(p)}p
                              \bigl(p^{-\kappa}-p^{-\sigma}\bigr).
    \]
   We note that the primes $p\le x$ contribute 
    \[
    2(\sigma-\kappa)\sum_{p\le x}\frac{\log p}{p} \ll 1,
    \]
    whereas the primes $p>e^{1/\kappa}$ also contribute $$\sum_{p \ge  e^{1/\kappa}} p^{-1-\kappa} \ll 1.$$ The conclusion then follows by restricting our sum to $p\in [x,e^{1/\kappa}].$
\end{proof}

\subsection{An explicit cutoff}
Our goal now is to recall a standard explicit cut off, which we shall use to perform explicit estimates. To this end, for $Y>1$ and $x>0$, we introduce
\begin{equation}\label{eq:averaged-cutoff}
 W_Y(x)=\int_{[0,1]^6}
 \1_{\{x\le Y^{1+(u_1+\cdots+u_6)/6}\}}
 \,du_1\cdots du_6.
\end{equation}
We recall that $0\le W_Y\le1$, with $W_Y(x)=1$ for $x\le Y$ and
$W_Y(x)=0$ for $x\ge Y^2$. We further have that $W_Y$ is continuously differentiable
with the decay estimate
\begin{equation}\label{eq:averaged-cutoff-derivative}
 |xW_Y'(x)|\le\frac{100}{\log Y}.
\end{equation}
Clearly, the sixth variable convolution is merely a technicality and shall be used only to ensure enough decay of the corresponding derivatives.
For $z\ne0$, we put
\begin{equation}\label{eq:perron-kernel}
 \Omega_Y(z)=\frac{Y^z}{z}
       \left(\frac{Y^{z/6}-1}{z\log Y/6}\right)^6,
\end{equation}
and note that this function has a single pole, at zero, of residue $1$ and furthermore 
for every $n$ and $c>0$,
\begin{equation}\label{eq:averaged-perron}
 W_Y(n)=\frac1{2\pi i}\int_{(c)}\Omega_Y(z)n^{-z}\,dz.
\end{equation}
The latter is standard and follows from the application of Perron's formula,
see
\cite[p.~1484, proof of Lemma~4.2 and (4.1)]{SoundWeakSubconvexity},
or \cite[Lemma~2.1 and its proof]{LamzouriClassNumbers}, with
$N=6$, $\lambda=(\log Y)/6$, and $y=Y/x$. Now summing over $n$ and picking $c = 1/\log Y$ and further making the change of variable $s = \kappa + z,$ we then get
\begin{equation}\label{eq:perron-weighted-sum}
 \sum_n\frac{g(n)W_Y(n)}{n^{1+\kappa}}
 =\frac1{2\pi i}\int_{(\kappa+1/\log Y)}
                 F(s)\Omega_Y(s-\kappa)\,ds.
\end{equation}
Moreover, we have from (\ref{eq:perron-kernel}) for all $v,$
\begin{equation}\label{eq:perron-positive-line}
 \left|\Omega_Y\left(\frac1{\log Y}+iv\right)\right|
 \ll\frac{\log Y}{(1+|v|\log Y)^7}.
\end{equation}

We are now ready to introduce the explicit cutoff in order to prove Proposition~\ref{prop:direct-elliott}, 
\begin{equation}\label{eq:perron-band}
    \Phi(x) = W_{X^{1/2}}(x) - W_{P^3}(x),
\end{equation}
which verifies the hypotheses of Proposition~\ref{elliott} by (\ref{eq:averaged-cutoff-derivative}). The next lemma lets us show the cutoff $\Phi$ picks up most of the weight in $F(\kappa)$ when $\kappa$ is roughly between $1/\log X$ and $1/\log P.$

\begin{lemma}\label{lem:perron-smoothing}
    We have for $Y \ge 3, 0 < \kappa < 1$ and $K := \kappa \log Y,$
    \begin{equation}\label{kappa-low}
        \begin{gathered}
            \sum_n \frac{g(n)W_Y(n)}{n^{1+\kappa}} \ll |F(\kappa)| K e^{\D^2(g, 1, Y, Y^{1/K})} \quad (K\le 1),
        \end{gathered}
    \end{equation}
    and 
    \begin{equation}\label{kappa-high}
        \begin{gathered}
            \sum_n \frac{g(n)W_Y(n)}{n^{1+\kappa}} = F(\kappa) \Big(1 + O(e^{-K/2}K^{-6})\Big) \quad (K\ge 1).
        \end{gathered}
    \end{equation}
    %In particular for $\kappa = (\log P\log X)^{-1/2},$ we have
    %\begin{equation}
    %    \sum_n \frac{g(n)\Phi(n)}{n^{1+\kappa}} = F(\kappa)\Big(1 + O()\Big)
    %\end{equation}
\end{lemma}
\begin{proof}
    We first prove \eqref{kappa-low}. It follows from \eqref{eq:perron-weighted-sum}, Lemma~\ref{lem:vertical-euler} and \eqref{eq:perron-positive-line} that
    \begin{equation}\label{eq:perron-cutoff-upper}
        \left|\sum_n\frac{g(n)W_Y(n)}{n^{1+\kappa}}\right|
        \ll\left|F\left(\kappa+\frac1{\log Y}\right)\right|.
    \end{equation}
    Combining with Lemma~\ref{lem:horizontal-G} gives \eqref{kappa-low}.
    \\ We now deal with \eqref{kappa-high}. We come back to \eqref{eq:perron-weighted-sum} and shift the
    contour to $\re s=\kappa/2$. Only the pole at $s=\kappa$ is
    crossed, with residue $F(\kappa)$ and the horizontal integrals vanish 
    by the decay of the kernel and the bound
    $|F(s)|\le\zeta(1+\kappa/2)$ in the strip. Thus
    \[
     \sum_n\frac{g(n)W_Y(n)}{n^{1+\kappa}}-F(\kappa)
     =\frac1{2\pi i}\int_{(\kappa/2)}
                      F(s)\Omega_Y(s-\kappa)\,ds.
    \]
    By Lemma~\ref{lem:vertical-euler} we have
    \begin{align*}
     \left|\Omega_Y\left(-\frac\kappa2+iv\right)\right|
     &\ll\frac{e^{-K/2}\log Y}{(K+|v|\log Y)^7},\\
     |F(\kappa/2+iv)|
     &\ll\left(1+\frac{|v|\log Y}{K}\right)^2|F(\kappa)|.
    \end{align*}
    Consequently, the contour integral is
    $O(e^{-K/2}K^{-6}|F(\kappa)|)$, which proves \eqref{kappa-high}.
\end{proof}

\subsection{A refined Elliott inequality}
We now estimate the error term in \eqref{multiplicative-elliott} in order to prove Proposition~\ref{prop:direct-elliott}.

\begin{lemma}\label{lem:direct-cutoff}
    Let $P$ be sufficiently large, let $X\ge P^{12}$, let $0<\kappa\le1/(3\log P)$ and $\mathcal D\ge0$, and assume that $\D^2(g, 1, P^3,  e^{1/\kappa}) \le\mathcal D.$ Uniformly for primes $p\le P$,
    \begin{equation}\label{eq:direct-dilation-bound}
     \begin{split}
     &\left|\sum_n\frac{g(n)}{n^{1+\kappa}}
                \bigl(\Phi(n)-p^{-\kappa}\Phi(pn)\bigr)\right|\\
     &\quad\ll  e^{\mathcal D}|G(\kappa)|\log p
     \left(\frac1{\log P}+
     \sum_q\frac{|1-g(q)|}{q^{1+\kappa}}
                   \bigl(1-q^{-1/(3\log P)}\bigr)\right),
     \end{split}
    \end{equation}
    with an absolute implied constant and  $\Phi$ is defined accordingly by \eqref{eq:perron-band}.
\end{lemma}

\begin{proof}
    For $Y\in\{P^3,X^{1/2}\}$ we put
    $\sigma=\kappa+1/\log Y$ and
using Perron's formula, we have
    \begin{equation}\label{eq:perron-dilation-identity}
     \begin{split}
     &\sum_n\frac{g(n)}{n^{1+\kappa}}
            \bigl(W_Y(n)-p^{-\kappa}W_Y(pn)\bigr)\\
     &\qquad=\frac1{2\pi i}\int_{(\sigma)}
            \frac{1-p^{-s}}s\,G(s)\Omega_Y(s-\kappa)\,ds.
     \end{split}
    \end{equation}
   We write $G(s)=G(\kappa)+G(s)-G(\kappa)$ and note $(1-p^{-s})/s=\int_1^p v^{-1-s}\,dv,$ so the contribution of
    $G(\kappa)$ in place of $G(s)$ is
    \[
     G(\kappa)\int_1^p v^{-1-\kappa}W_Y(v)\,dv.
    \]
    This follows by interchanging absolutely convergent integrals and
    applying \eqref{eq:averaged-perron}. On subtracting the two cutoff
    expressions the constant (pole) contributions crucially cancel since
    \begin{equation}\label{eq:direct-pole-cancellation}
     G(\kappa)\int_1^p\frac{\Phi(v)}{v^{1+\kappa}}\,dv=0,
     \qquad p\le P,
    \end{equation}
    because $\Phi$ vanishes on $[1,P^3].$ We may therefore replace
    $G(s)$ by $G(s)-G(\kappa)$ in both contour integrals, provided we give a satisfactory pointwise estimate for $|G(s)-G(\kappa)|.$
    To this end, we first note that for all $v\in\R,$ by Lemma~\ref{lem:vertical-euler} we have %Lemmas~\ref{lem:vertical-euler, lem:horizontal-G} gives
    \begin{equation}\label{eq:G-flat}
     |G(\sigma+iv)|\ll(1+|v|\log Y)^3|G(\sigma)|\ll  e^{\mathcal D}(1+|v|\log Y)^3|G(\kappa)|.
    \end{equation}
     To see that the second inequality holds, we note that if $\kappa\le1/\log Y$ (which is the case for $Y=P^3,$ say), we use can Lemma~\ref{lem:horizontal-G} together with the bound $\D^2(g,1,Y, e^{1/\kappa})\le\D^2(g, 1, P^3,  e^{1/\kappa})\le\mathcal D$. Now if $\kappa>1/\log Y$, then $\kappa<\sigma\le2\kappa$, and directly the Euler product gives
    \[
     \log\frac{|G(\sigma)|}{|G(\kappa)|}
     \le\log\frac\sigma\kappa+\sum_p\bigl(p^{-1-\kappa}-p^{-1-\sigma}\bigr)+O(1)
     =O(1).
    \]
    
    We now show that, for every $v\in\R$,
    \begin{equation}\label{eq:finite-difference-G}
     \begin{split}
     |G(\sigma+iv)-G(\kappa)|
     &\ll  e^{\mathcal D}|G(\kappa)|(1+|v|\log Y)^4\\
     &\quad\times\left(\frac1{\log Y}+
         \sum_p\frac{|1-g(p)|}{p^{1+\kappa}}
                             \bigl(1-p^{-1/\log Y}\bigr)\right).
     \end{split}
    \end{equation}
    To this end, we fix an absolute $c_0\in(0,1)$ sufficiently small that $s\zeta(1+s)$
    is holomorphic and non-zero for $|s|\le3c_0$.
    If $|v| \ge c_0,$
    the estimate follows directly from \eqref{eq:G-flat}, since
    $(1+|v|\log Y)/\log Y>c_0$.
    We therefore suppose that $|v|\le c_0$ and assume without loss of generality that $\kappa+1/\log Y\le c_0.$ Since the logarithmic derivative of
    $s\zeta(1+s)$ is bounded, we have 
    \[
     \begin{split}
     |\log G(\sigma+iv)-\log G(\kappa)|
     \ll{}&|\sigma+iv-\kappa|\\
     &+\sum_p\sum_{j\ge1}\frac{|1-g(p)^j|}{j p^{j(1+\kappa)}}
                 \left|1-p^{-j(1/\log Y+iv)}\right|.
     \end{split}
    \]
    For each $p,j$, we have
    \[
     \left|1-p^{-j(1/\log Y+iv)}\right|
     \le(1+|v|\log Y)(1-p^{-j/\log Y})
     \le j(1+|v|\log Y)(1-p^{-1/\log Y}).
    \]
    where the first inequality follows from $u e^{-u}\le1-e^{-u}$ for
    $u\ge0$. Using further $|1-g(p)^j|\le j|1-g(p)|$ and summing the geometric
    series, we obtain
    \[
     \begin{split}
     |\log G(\sigma+iv)-\log G(\kappa)|
     &\ll (1+|v|\log Y)\\
     &\quad\times\left(\frac1{\log Y}+
       \sum_p\frac{|1-g(p)|}{p^{1+\kappa}}
                        (1-p^{-1/\log Y})\right).
     \end{split}
    \]
    Now using $|e^z-e^w|\le\max(e^{\re z},e^{\re w})|z-w|$, and combining this
    with \eqref{eq:G-flat} we establish \eqref{eq:finite-difference-G}.

    Applying \eqref{eq:finite-difference-G},
    $|(1-p^{-s})/s|\le\log p$, and
    \eqref{eq:perron-positive-line} and noting that the fourth power in
    \eqref{eq:finite-difference-G} is absorbed by the seventh-power
 decay of the kernel, we have
    \[
     \int_\R\frac{\log Y}{(1+|v|\log Y)^3}\,dv\ll1.
    \]
    Finally since $\log Y\ge3\log P$, so
    $1-p^{-1/\log Y}\le1-p^{-1/(3\log P)}$ and 
    summing the bounds for the two cutoffs delivers 
    \eqref{eq:direct-dilation-bound}.
\end{proof}

\subsection{Completion of the proofs}
\begin{proof}[Proof of Proposition \ref{prop:direct-elliott}]
   We are now ready to complete the proofs of our main results. Let $\kappa:=1/(A\log P)$ and let $\Phi$ be given by \eqref{eq:perron-band}. We take large $B\ge4$ and $C\ge3$, so that $AB\ge12$ and $X^{1/2}\ge P^6.$ Now $\Phi$ is supported on $[P^3,X]$ and clearly $0\le\Phi\le1$ and by \eqref{eq:averaged-cutoff-derivative} we also have $|\Phi'(u)|\le4/(u\log P)\le1/u$. Hence Proposition~\ref{elliott} applies, and combining \eqref{elliott-eq} with \eqref{multiplicative-elliott} we get
    \begin{equation}\label{start-6.0}
        \left|\sum_n\frac{g(n)\Phi(n)}{n^{1+\kappa}}\right|
        \sum_{\substack{p\le P\\p\nmid m}} |1-g(p)| \ll 
        \frac{(1+|t|)\Delta}{\kappa\log P} + 
        \sum_{\substack{p\le P\\p\nmid m}}
        \left|\sum_n\frac{g(n)}{n^{1+\kappa}} \big(\Phi(n) - p^{-\kappa} \Phi(pn)\big)\right|.
    \end{equation}
    We now apply Lemma~\ref{lem:perron-smoothing} repeatedly twice. For $Y=X^{1/2}$ we have $K=\kappa\log X^{1/2}=B/2\ge1$, while for $Y=P^3$ we have $K=3/A\le1$ and $Y^{1/K}=P^A= e^{1/\kappa}$. Since $$\D^2(g, 1, P^3,  e^{1/\kappa})=\D(f,\chi(n)n^{it};P^3,P^A)^2\le\mathcal D,$$ we obtain
    \begin{equation}\label{estim-smoothed-sum}
    \begin{split}
        \sum_n\frac{g(n)\Phi(n)}{n^{1+\kappa}} &= \sum_n\frac{g(n)W_{X^{1/2}}(n)}{n^{1+\kappa}} - \sum_n\frac{g(n)W_{P^3}(n)}{n^{1+\kappa}}
        \\&= F(\kappa) \bra{1 + O\bra{ e^{-B/4}B^{-6}} + O\bra{A^{-1} e^{\mathcal D}}}.
    \end{split}
    \end{equation}
    Choosing first $B$ and then $C$ sufficiently large (absolute), the assumption $A\ge C e^{\mathcal D}$ ensures that
    \begin{equation}
        \left|\sum_n\frac{g(n)\Phi(n)}{n^{1+\kappa}}\right| \ge \frac12|F(\kappa)|.
    \end{equation}
    Next, since $\kappa\le1/(3\log P)$ and $X\ge P^{12}$, Lemma~\ref{lem:direct-cutoff} applies to every $p\le P$ and together with $\sum_{p\le P}\log p\le2P$ and \eqref{start-6.0}, this yields
    \begin{equation}
     \begin{split}
         \sum_{\substack{p\le P\\p\nmid m}} |1-g(p)| \ll{}& \frac{(1+|t|)\Delta}{|F(\kappa)|\,\kappa \log P}\\
         &+ \frac{ e^{\mathcal D}|G(\kappa)| P}{|F(\kappa)|} \left(\frac1{\log P}+
            \sum_p\frac{|1-g(p)|}{p^{1+\kappa}}
                   \bigl(1-p^{-1/(3\log P)}\bigr)\right).
     \end{split}
    \end{equation}
    Using the bound $|F(\kappa)| \gg \kappa$ from Lemma~\ref{lem:vertical-euler} and the definitions $G(\kappa) = \kappa F(\kappa)$ and $\kappa = 1/(A\log P)$ gives \eqref{eq:direct-elliott-bound}, after enlarging $C$ if necessary, since the primes $p\mid m$ contribute at most $\omega(m)$ to the left-hand side.
\end{proof}

\begin{proof}[Proof of Theorem \ref{thm:subpolynomial}]
    Let $B,C$ be the absolute constants of Proposition~\ref{prop:direct-elliott}, and let $\eps_1=\eps_1(\delta)$ and $\mathcal D=\mathcal D(\delta)$ be as in Corollary~\ref{cor:stabilization-fixed}. Let $A=A(\delta)$ be such that  $A\ge C e^{\mathcal D}$ and
    \begin{equation}\label{eq:theta-choice}
        \theta:=\frac{C e^{\mathcal D}\log(1+A/3)}{A}\le\frac12.
    \end{equation}
    We now take
    \[
        \eps_0(\delta):=\min\left\{\eps_1(\delta),\,\frac{1}{2AB}\right\}.
    \]
   We note that these quantities depend only on $\delta$. We now assume that \eqref{eq:l2-density} and \eqref{eq:subpoly-hyp} hold. Then $\mathcal N_f(x)\le x^{\eps_0+\eta(x)}$ for all large $x$, where $\eta(x)\ge0$ and $\eta(x)\to0$. By Corollary~\ref{cor:stabilization-fixed} there are a primitive Dirichlet character $\chi$, of conductor $m$, and $t\in\R$, which we fix from now on, such that
    \begin{equation}\label{eq:stab-fixed-applied}
        \D(f,\chi(n)n^{it};x^{3\eps_0},x)^2\le\mathcal D
    \end{equation}
    for all large $x$. We let $g(n) = f(n)\conj{\chi(n)}n^{-it}$ and introduce the quantities
    \begin{equation}
        \begin{gathered}
            V(P) := \sum_{p\le P} \abs{f(p) - \chi(p) p^{it}}, \qquad 
            S(P) := \sum_{p\le P} \frac{\abs{f(p) - \chi(p) p^{it}}}{p}.
        \end{gathered}
    \end{equation}
    Note that $|f(p)-\chi(p)p^{it}|\le|1-g(p)|$ for every prime $p$, with equality if $p\nmid m$.

    Let $P$ be large and set $X:=P^{AB}$. Applying now \eqref{eq:stab-fixed-applied} with $x=P^A$, and noting that $P^{3A\eps_0}\le P^3$, we get $\D(f,\chi(n)n^{it};P^3,P^A)^2\le\mathcal D$. Furthermore, for $\eta^*(P):=\sup_{u\ge P}\eta(u),$ using the trivial bound $\mathcal N_f(u)\le\sqrt u$ for $u\le P$ and $AB\eps_0\le1/2$, we have
    \[
        \mathcal N_f(u)\le\Delta:=P^{1/2+AB\eta^*(P)}\qquad(1\le u\le X),
    \]
    and $\Delta\le P^{2/3}$ for all large $P$. Proposition~\ref{prop:direct-elliott} therefore gives, for all sufficiently large $P$,
    \begin{equation}\label{eq:V-contracting}
     \begin{split}
     V(P)\le{}&\frac{C e^{\mathcal D}P}{A\log P}
        \sum_p\frac{|1-g(p)|}{p^{1+1/(A\log P)}}
                           \bigl(1-p^{-1/(3\log P)}\bigr)\\
     &+\frac{C e^{\mathcal D}P}{A(\log P)^2}+O_f\bigl(P^{3/4}\bigr).
     \end{split}
    \end{equation}
    Now, by partial summation we have for all $T,$
    \[
     S(T)=\frac{V(T)}T+\int_2^T\frac{V(P)}{P^2}\,dP.
    \]
    Inserting \eqref{eq:V-contracting} for large $P$, and the trivial bound $V(P)\le2\pi(P)$ for the remaining $P$, yields
    \begin{equation}
     \begin{split}
     S(T)\le C_{f}+\frac{C e^{\mathcal D}}{A}\sum_p\frac{|1-g(p)|}{p}
     \int_{2}^T p^{-1/(A\log P)}
           \bigl(1-p^{-1/(3\log P)}\bigr)\frac{dP}{P\log P}.
     \end{split}
    \end{equation}
    For each prime $p$, we make a substitution $u=\frac{\log p}{A \log P}$, so that 
    \begin{equation}\label{eq:integrated-prime-sum}
     \begin{split}
     S(T)\le C_{f}+\frac{C e^{\mathcal D}}{A}\sum_p\frac{|1-g(p)|}{p}
     \int_{\frac{\log p}{A \log T}}^\infty e^{-u}
           \bigl(1-e^{-uA/3}\bigr)\frac{du}{u}.
     \end{split}
    \end{equation}
    For $p\ge T^A,$ the inner integral is $\le p^{-1/A\log T},$ and the contribution of terms $p \ge T^A$ to \eqref{eq:integrated-prime-sum} is therefore at most
    \begin{equation}
        \ll \frac{C e^{\mathcal D}}{A}\sum_{p\ge T^A} p^{-1 - 1/A\log T} \ll 1,
    \end{equation}
    since $A\ge C e^{\mathcal D}$. For $p< T^A,$ we bound the integral by $\int_0^\infty e^{-u}(1-e^{-uA/3})\,\frac{du}u=\log(1+A/3)$, and we recall that $|1-g(p)|=|f(p)-\chi(p)p^{it}|$ for $p\nmid m$. Thus \eqref{eq:integrated-prime-sum} becomes
    \begin{equation}\label{eq:contraction-target}
     S(T)\le C_{f}+\theta\,S(T^A),
    \end{equation}
    where $C_f$ also absorbs $\theta\sum_{q\mid m}1/q$, and $\theta\le1/2$ by \eqref{eq:theta-choice}. Iterating the latter gives
    \[
     S(T)\le C_{f}\sum_{j=0}^{k-1}\theta^j
                       +\theta^kS(T^{A^k}).
    \]
    Since $|f(p)-\chi(p)p^{it}|\le2$, Mertens' theorem yeilds the inequality
    $S(T^{A^k})\ll1+\log\log T+k\log A$. Letting $k\to\infty$, we
    conclude that $S(T)\le2C_f$ for all $T$. This establishes
    \eqref{eq:subpoly-conclusion}.
\end{proof}

\section{Sharpness of the results}\label{sec:example}
Our goal now is to exhibit an example satisfying \ref{thm:nonpretentious-example}, which shows that in Theorem~\ref{thm:subpolynomial} the power saving \eqref{eq:subpoly-hyp} has to hold for all sufficiently large scales and hence, the $\limsup$ there cannot be replaced by the $\liminf$.   
We begin with the following basic estimate handling completely multiplicative functions, which vary on finite set of primes.
\begin{lemma}\label{lem:finite-modifications}
Suppose $g:\N\to\{-1,1\}$ is completely multiplicative and $g(p)=1$
outside a finite set of $r$ primes. Then
\[
 \mathcal N_g(x)\ll(\log(2x))^{r+1}.
\]
\end{lemma}
\begin{proof}
Let $g=1*h$ and note that $h$ is supported on integers composed of exceptional primes.
The inner sum can be rewritten as
\[
 \sum_{d\le x}h(d)\sum_{m\le x/d}e(dm\alpha)
\]
and invariance of the $L_1$ norm under $\alpha\mapsto d\alpha$ gives
\[
 \mathcal N_g(x)\ll\log(2x)\sum_{d\le x}|h(d)|,
\]
where the term $O(\log(2x))$ comes from bounding geometric sum
by integrating
$|\sum_{m\le u}e(m\alpha)|\ll\min(u,\|\alpha\|^{-1})$.
Finally,
\[
 \sum_{d\le x}|h(d)|
 \le\prod_{p\ \mathrm{exceptional}}
       \left(1+2\left\lfloor\frac{\log x}{\log p}\right\rfloor\right)
 \ll(\log(2x))^r,
\]
and the result follows.
\end{proof}
The following classical bound is a direct consequence of the celebrated result of Vaughan \cite{Vaughan}.
\begin{lemma}\label{lem:prime-polynomial}
For every sufficiently large integer $N$,
\begin{equation}\label{eq:prime-polynomial-lower}
 \int_0^1\left|\sum_{N<p\le2N}e(p\alpha)\right|\,d\alpha
 \gg\frac{\sqrt {N}}{(\log N)^2}.
\end{equation}
\end{lemma}

\subsection{The construction}
\begin{proof}[Proof of Theorem~\ref{thm:nonpretentious-example}]
We construct the function by changing its values on widely separated
finite blocks of primes
\[
 \mathcal P_j=\{p:y_j<p\le z_j\},\qquad j\ge2,
\]
Start with $z_1=1$ and
suppose the preceding blocks have been chosen, and let $g_{j-1}$ be
the completely multiplicative function with $g_{j-1}(p)=-1$ on their
union and $g_{j-1}(p)=1$ elsewhere. Lemma~\ref{lem:finite-modifications}
allows us to choose $y_j>z_{j-1}$ sufficiently large that
\begin{equation}\label{eq:example-Y}
 \log\log y_j\ge j^2,
 \qquad
 \mathcal N_{g_{j-1}}(x)
 \le\min\{x^{1/j},\sqrt{x}/(\log x)^3\}
 \quad(x\ge y_j).
\end{equation}
We also impose the condition
\begin{equation}\label{eq:example-block-start}
 \sum_{y_j<p\le2y_j}\frac1p<\frac1j,
 \qquad y_j>j.
\end{equation}
The first condition holds for all sufficiently large $y_j.$ We choose $z_j$ to be the first
prime for which
\[
 \sum_{y_j<p\le z_j}\frac1p\ge\frac1j.
\]
Since the final summand is at most $1/y_j<1/j$, we have
\begin{equation}\label{eq:example-block-mass}
 \frac1j\le\sum_{p\in\mathcal P_j}\frac1p\le\frac2j.
\end{equation}
By \eqref{eq:example-block-start}, we clearly have $z_j>2y_j$ and finally
define $f$ at primes by
\[
 f(p)=\begin{cases}
 -1,&p\in\displaystyle\bigcup_{j\ge2}\mathcal P_j,\\
 1,&\text{otherwise},
 \end{cases}
\]
and extend completely multiplicatively.
\subsection*{Nonpretentiousness.} We  claim that the function $f$ from above is non-pretentious. 
 Indeed, we note that for $y_j\le x<y_{j+1}$,
\[
 \sum_{p\le x,\ p\in\bigcup_j\mathcal P_j}1/p\ll\sum_{k=2}^j\frac1k\ll\log j,
 \qquad \log\log x\ge j^2,
\]
and the left hand side tends to infinity with $x,$ so that 
$
 \D(f,1;\infty)=\infty,
$
and $\D(f,\chi_0;\infty)=\infty,$ where $\chi_0$ is a principle character.
For any fixed twist $\chi(n)n^{it}$, except a principal character
with $t=0$, the prime number theorem in
arithmetic progressions gives
\begin{equation}\label{eq:fixed-twist-distance}
 \D(1,\chi(n)n^{it};x)^2=\log\log x+O_{\chi,t}(1).
\end{equation}
Replacing $1$ by $f(p)$
on the selected primes changes the squared distance by at most
$o(\log\log x)$ and consequently $ \D(1,\chi(n)n^{it};x)^2=(1+o(1))\log\log x$ for every non-principle $\chi.$ We can now conclude that $f$ is (strongly) nonpretentious.

\subsection*{The bound at every scale.}
We now fix $y_j\le x<y_{j+1}$ and note that for $n\le x$, the values $f(n)$ and
$g_{j-1}(n)$ agree unless a prime from $\mathcal P_j$ divides $n$.
Consequently,
\[
 \sum_{n\le x}|f(n)-g_{j-1}(n)|^2
 \le4\sum_{p\in\mathcal P_j}\left\lfloor\frac xp\right\rfloor
 \le\frac{8x}{j}.
\]
Applying the triangle inequality, Parseval and Cauchy--Schwarz yield 
\begin{equation}\label{eq:example-all-scales}
 \mathcal N_f(x)
 \le\mathcal N_{g_{j-1}}(x)+\sqrt{8x/j}
 \le\frac{\sqrt{x}}{(\log x)^3}+\sqrt{\frac{8x}{j}}.
\end{equation}
Since $x\to\infty$ forces $j$ to tend to infinity, this proves
\eqref{eq:example-uniform} for all sufficiently large $x$.
At the point $x=y_j$, the current block has not yet contributed, and hence
\[
 \mathcal N_f(y_j)=\mathcal N_{g_{j-1}}(y_j)\le y_j^{1/j}.
\]
 Thus, the liminf bound in
\eqref{eq:example-exponents} follows by taking $j\to\infty.$

\subsection*{The lower bound along a subsequence.}
Since $z_j>2y_j$, every prime in $(y_j,2y_j]$ belongs to the current
block. None has a nontrivial multiple at most $2y_j$, and no later
block contributes in this range. Thus
\begin{equation}\label{eq:example-exact-spike}
 S_f(2y_j,\alpha)=S_{g_{j-1}}(2y_j,\alpha)
                    -2\sum_{y_j<p\le2y_j}e(p\alpha).
\end{equation}
Applying Lemma~\ref{lem:prime-polynomial}, the triangle
inequality, and \eqref{eq:example-Y}, we obtain
\[
 \mathcal N_f(2y_j)
 \ge2\left\|\sum_{y_j<p\le2y_j}e(p\alpha)\right\|_{L^1([0,1])}
           -\mathcal N_{g_{j-1}}(2y_j)
 \gg\frac{\sqrt{y_j}}{(\log y_j)^2}.
\]
Together with the trivial bound
$\mathcal N_f(x)\le\sqrt{\lfloor x\rfloor}$, this proves the limsup claim
in \eqref{eq:example-exponents}. Moreover, for every fixed $a<1/2$,
\[
 \frac{\mathcal N_f(2y_j)}{(2y_j)^a}\longrightarrow\infty,
\]
and the final conclusion follows.
\end{proof}

\section*{Acknowledgements}

The authors thank R\'egis de la Bret\`eche for his insightful comments on
earlier drafts of the paper.\\
The main results were established in January 2026 during the second
author's visit to IHP and have been presented during various seminars/workshops between June and September, 2026 (before the current AI ``revolution"). He thanks the institute for its warm hospitality. We are also grateful to the participants of the Heilbronn focused research workshop (in February, 2026), in particular, Benjamin Bedert and Thomas Bloom, for valuable feedback and discussions. 
Part of the work was carried out during the Universal Statistics in
Number Theory semester in Montreal. The authors are very grateful to the CRM and
Universit\'e de Montr\'eal for their hospitality during their stay. \\
The first author thanks the CRM as well as the FSMP for their financial support during her stay in Montreal.

The AI tools were used at the final stage to track down, check and correct numerous explicit dependencies of the parameters. The ideas and calculations are due to the authors.

\end{document}